\documentclass{amsart}
\usepackage{amsmath,amsthm,amssymb,amsfonts}
\ifx\pdfoutput\undefined \usepackage[hypertex]{hyperref} \else \usepackage[pdftex,pdfstartview=FitH,pdfpagemode=none]{hyperref} \fi
\numberwithin{equation}{section}
\newcommand{\IGNORE}[1]{}
\newcommand{\tr}{\operatorname{tr}}
\newcommand{\haarG}{\mu_{\mathrm{Haar}}}

\newcommand{\z}{\mathfrak{Z}}

\renewcommand{\t}{\mathfrak{t}}
\newcommand{\T}{\mathbf{T}}

\newcommand{\Ad}{\operatorname{Ad}}

\newcommand{\disc}{\operatorname{disc}}
\newcommand{\h}{\mathfrak{h}}

\newcommand{\D}{\mathbf{D}}

\newcommand{\g}{\mathfrak{g}}
\newcommand{\C}{\mathbb{C}}

\newcommand{\vol}{\mathrm{vol}}

\newcommand{\G}{\mathbf{G}}

\newcommand{\GL}{\operatorname{GL}}

\newcommand{\PGL}{\operatorname{PGL}}

\newcommand{\nope}{}

\newcommand{\SL}{\operatorname{SL}}
\newcommand{\diag}{\operatorname{diag}}

\newcommand{\order}{\mathscr{O}}
\newcommand{\N}{\mathbf{N}}
\newcommand{\R}{\mathbb{R}}
\renewcommand{\H}{\mathbb{H}}

\newcommand{\Z}{\mathbb{Z}}
\newcommand{\Q}{\mathbb{Q}}

\newcommand{\adele}{\mathbb{A}}

\newcommand{\Lie}{\operatorname{Lie}}
\DeclareFontFamily{OT1}{rsfs}{}
\DeclareFontShape{OT1}{rsfs}{n}{it}{<-> rsfs10}{}
\DeclareMathAlphabet{\mathscr}{OT1}{rsfs}{n}{it}

\newtheorem{Theorem}{Theorem}[section]
\newtheorem*{Theorem*}{Theorem}
\newtheorem{Proposition}[Theorem]{Proposition}
\newtheorem{Lemma}[Theorem]{Lemma}
\newtheorem{Corollary}[Theorem]{Corollary}
\newtheorem{Conjecture}[Theorem]{Conjecture}
\theoremstyle{definition}
\newtheorem{Definition}[Theorem]{Definition}
\newtheorem{Example}[Theorem]{Example}

\newtheorem{Remark}[Theorem]{Remark}
\newtheorem*{Acknowledgements}{Acknowledgements}

\newcommand\equ{\eqref}
\newcommand {\norm}[1] {\left\| {#1} \right\|}
\newcommand{\weaklygoesto}{\stackrel{w*}{\goesto}}
\newcommand {\goesto} {\longrightarrow}
\newcommand {\absolute}[1] {\left| {#1} \right|}
\newcommand {\floorof}[1] {\left\lfloor {#1} \right\rfloor}
\newcommand\htop{\operatorname{h_{top}}}

\begin{document}

\title[Distribution of periodic torus orbits]{The distribution of periodic torus orbits on homogeneous spaces}
\author[M.~Einsiedler, E.~Lindenstrauss, Ph.~Michel, A.~Venkatesh]{Manfred Einsiedler, Elon Lindenstrauss, Philippe Michel and Akshay Venkatesh}
\thanks{This research has been supported by the Clay Mathematics Institute (M.E. by a Clay Research Scholarship, E.L. and A.V. by fellowships),
 by the ``R\'eseau arithm\'etique des Pyr\'en\'ees'' of R\'egion Languedoc-Roussillon
and the  European TMR Network ``Arithmetic Algebraic Geometry'' (Ph.M.); and by the NSF (DMS grants 0509350 (M.E.), 0500205 (E.L.), 02045606 (A.V.)
 and an FRG collaborative grant.)}

\date{7/27/2006}

\begin{abstract}
We prove results towards the equidistribution of certain families of periodic
torus orbits on homogeneous spaces, with particular focus
on the case of the diagonal torus acting on quotients of $\PGL_n(\R)$.
After attaching to each periodic orbit an integral invariant (the discriminant)
our results have the following flavour:
certain standard conjectures about the distribution of such orbits
hold up to exceptional sets of at most $O(\Delta ^{\epsilon})$
orbits of discriminant $\leq \Delta$. The proof relies on the well-separatedness of periodic orbits together with measure rigidity for torus actions.
We also give examples of sequences of periodic orbits of this action that
 fail to become equidistributed, even in higher rank.

We give an application of our results to sharpen a theorem of Minkowski on ideal classes in totally real number fields of cubic and higher degrees.
\end{abstract}
\maketitle

\tableofcontents
\newcounter{savedenumi}

\section{Introduction.} \label{RealIntro}

\subsection{Periodic torus orbits}\label{opening subsection}
Let $G = \mathbf G (\R)$ be a real algebraic group defined over $\Q$ (later we will specialize to the case of $\R$-split groups), and $\Gamma < G$ an arithmetic lattice.\footnote {I.e. a lattice in $G$ contained in $\mathbf G (\Q)$.}
Let $L$ be a closed subgroup of $G$.
We say that an $L$-orbit on $\Gamma \backslash G$ is \emph{periodic} if it possesses a finite $L$-invariant measure\footnote {This is not a completely standard definition, but is quite natural.}.
Periodic orbits of subgroups $L<G$ on such arithmetic quotients
have been studied in various context by several authors, and have proved to be a fruitful meeting ground of dynamics and number theory.
We note that there is nothing special about \emph{real} algebraic groups: we could have equally well considered $p$-adic or $S$-adic\footnote {I.e. products of real and $p$-adic groups (possibly for several different primes $p$).} groups. We focus on real groups only to simplify matters and for concreteness of exposition.

The multiple facets of the study of periodic orbits on arithmetic quotients can already be seen in the simple but important case of
$\G = \PGL_2$.
In this case, it is well-known that the following objects are essentially in bijection:
\begin{enumerate}
\item Closed geodesics on the {\em modular surface} $ \SL_2(\Z) \backslash \mathbb{H}$.
\item Periodic orbits of the diagonal subgroup of $\PGL_2(\R)$ on the arithmetic quotient $\PGL_2(\Z) \backslash \PGL_2(\R)$.
\item \label{three} Ideal classes in real quadratic orders (e.g. $\Z[\sqrt{d}]$ for $d > 0$).
\end{enumerate}

A natural question one is led to, when studying periodic $L$-orbits on arithmetic quotients, is the following:

\medskip

\emph {\textbf{Basic Question:} to what extent do larger and larger periodic $L$-orbits fill out more and more of $\Gamma \backslash G$?}

\medskip

If the group $L$ is generated by its \emph{unipotent}{}\footnote {An element $g$ of a (linear) algebraic group $\G(K)$ is said to be unipotent if for some faithful representation $ \rho$ of $\G$, all eigenvalues of $\rho(g)$ over the algebraic closure of $K$ are equal to $1$.} elements, Ratner's measure classification theorem \cite{Ratner-Annals}, together with a ``linearization'' technique developed in the work of Ratner \cite{Ratner-Duke}, Dani-Margulis \cite{DM}, Shah \cite{Shah-uniformly-distributed} and elsewhere, can be used to give an essentially complete answer to the basic question above: larger and larger periodic $L$-orbits become equidistributed in $\Gamma \backslash G$
unless there are certain obvious obstacles present (see e.g.
\cite{Eskin-Oh} and \cite{Mozes-Shah} for such results with arithmetic consequences).

Our aim in this paper and its sequels (currently two are planned) is to study periodic orbits of a maximal $\R$-split torus $H<G$ on $G/\Gamma$ for $G= \mathbf G (\R)$ an $ \R$-split group. Such groups $H$ have by definition no unipotent elements whatsoever, and their action is much less well-understood\footnote {At least for $\dim H \geq 2$.
}.
While much of our discussion is quite general, we focus in the introduction and the later sections of this paper on the following two concrete examples:
\begin{enumerate}
\item [(L-1)] $\G = \PGL_n$, $G = \PGL_n(\R)$, $\Gamma = \PGL_n(\Z)$
\item [(L-2)] $D_{\Q}$ a degree $n$ division algebra over $\Q$, with $D_{\Q} \otimes_{\Q} \R \cong M_n(\R)$, $\G$ the group associated to invertible elements of $D$ modulo center,
$G = \G(\R) \cong \PGL_n(\R)$ and $\Gamma < G$ the lattice associated with an order $\mathcal{O} _ D$ in $D_{\Q}$.\footnote{See \S\ref{Section algebras} for a more complete discussion.}
\end{enumerate}
\noindent Identifying $G \cong \PGL _ n (\R)$, in both cases we can
take $H$ to be the group of $n \times n$ diagonal matrices\footnote
{Considered as elements in $\PGL _ n (\R)$, i.e. with proportional
matrices identified.}. Both examples come from a common family of
lattices in $\PGL_n(\R)$ arising from central simple algebras; in
the first case $\Gamma \backslash G$ is not compact, in the second
it is.

\subsection{Discriminant, shape and volume of periodic orbits.} \label{subsec:DSV}
Let $H$ be an $\R$-split maximal torus in $G$ as above\footnote {A maximal torus in $G$ is also called a Cartan subgroup, and we shall use the two terms interchangeably.}. Note that since $H$ is abelian, an orbit $xH = \Gamma g H$ of $H$ is periodic if and only if $xH$ is a compact subset of $\Gamma \backslash G$.

To such a periodic $H$-orbit we can attach several invariants. For example, we can look at the ``shape'' of the orbit, i.e. the stabilizer $g ^{-1} \Gamma g \cap H$ of $x$, which is a lattice
in  $H$. In particular, fixing a Haar measure on $H$, we can consider the \emph{volume} (or regulator) of the periodic orbit $xH$,
which by definition is the covolume  of $g ^{-1} \Gamma g \cap H$ in $H$.
We obtain another less obvious but nontheless important invariant of
the periodic orbit, the \emph{discriminant} of $xH$ (denoted $\disc
(xH)$), by associating to each periodic orbit a rational point on an
appropriate variety and looking at the denominator of this point. By
definition the discriminant is always a positive integer.

It turns out that the relative sizes of the discriminant and regulator of an orbits play a crucial role in the study of these orbits. In general, discriminant and regulator satisfy the following relations (see Proposition~\ref{proposition about regulator} below): there is some $c$ (which can be taken to be $1/2+ \epsilon$ for the arithmetic quotients given in (L-1) and (L-2) above) so that for any periodic $H$-orbit $xH$,
\begin{equation} \label{introduction regulator bound}
\log \disc (x H) \nope \ll \vol (x H) \nope \ll \disc (x H) ^c
\end{equation}
Our results, as well as the examples we give in Section \ref{examples section}, all give credence to the general principle that the {\em bigger the volume of a periodic orbit relative to its discriminant, the better the orbit is behaved.}

By considering periodic orbits individually, as we do in this paper, one does not take into account an  important hidden symmetry of the problem. Returning to the example of periodic $H$-orbits
on $\PGL_2(\Z) \backslash \PGL_2(\R)$, the ``volume'' of a periodic orbit is simply its length, or equivalently the length of the corresponding closed geodesics on $\SL (2, \Z) \backslash \H$, and the discriminant coincides with the discriminant of the associated real quadratic other.
As is well-known, the length spectrum of closed geodesics in this case is far from simple.

A similar property holds for all $n$, both for $ \Gamma =
\PGL_n(\Z)$ and more general lattices constructed from central
simple algebras and division algebras. The periodic $H$ orbits
naturally come in packets, with all orbits in a packet sharing the
same discriminant, regulator, and even shape. These packets can be
understood as projections to $\Gamma \backslash G$ of orbits of
adelic $\Q$-tori on $\G(\Q) \backslash \G(\adele)$. The compact
orbits belonging to a single packet are therefore parameterized by a
finite abelian group -- a suitable class group.

One important property of these packets is that their total volume (which is simply the volume of any single orbit in the packet times the number of such periodic orbits) is equal to $D ^{1/2 + o(1)}$ where $D$ is the discriminant of the packet --- i.e.  essentially equal to
the upper bound in \equ{introduction regulator bound}. These packets are considered in detail in \cite{ELMV2, ELMV3}.
A coarser grouping of periodic orbits can be obtained by using Hecke correspondences, and this is analyzed in \cite{Benoist-Oh}.

\subsection{Conjectures on higher rank rigidity.}\label{sec:higher rank rigidity}
The dynamics of the group $H$ is drastically different in the rank one and higher rank cases. In the rank one case, for example $\PGL_2( \Z) \backslash \PGL_2(\R)$,
the closing lemma assures us that each periodic orbit individually can be distributed
in an almost arbitrary way; the situation is very different, however, if one considers packets of periodic orbits as will be explained below.

In the higher rank cases, the dynamics is much more rigid. For example, G. A. Margulis has highlighted the following conjecture,
which is equivalent to a conjecture about the values of products of linear forms formulated by Cassels and Swinerton-Dyer in 1955 \cite{Cassels-Swinnerton-Dyer}:

\begin{Conjecture} \label{bounded orbits conjecture}
Let $H$ be the subgroup of diagonal matrices in $\PGL_n(\R)$ for $n \geq 3$. Any bounded\footnote{I.e. an orbit with compact closure.} $H$-orbit in $\PGL_n(\Z) \backslash \PGL_n(\R)$ is closed.
\end{Conjecture}

There are also related conjectures regarding invariant measures, due to H. Furstenberg (unpublished), A. Katok and R. Spatzier \cite{Katok-Spatzier}, and G. A. Margulis \cite{Margulis-conjectures}; a fairly general variant of these conjectures can be found in \cite[Conjecture 2.4]{ Einsiedler-Lindenstrauss-ICM}. Concretely one expects the following:

\begin{Conjecture} \label{measure classification conjecture}
Let $n \geq 3$, $G \cong \PGL_n(\R)$, $\Gamma < G$ a lattice as in
examples  (L-1) or (L-2) above\footnote {More generally, any lattice
coming from a central simple algebra over $\Q$; examples due to Rees
\cite{Rees-example} (a more accessible source is \cite [Section 9]{
Einsiedler-Katok}) show that in order to cover more general (even
cocompact) lattices in $G$ the conjecture should be reformulated to
allow for rank one factors (see e.g. \cite{Margulis-conjectures} or
\cite{Einsiedler-Lindenstrauss-ICM}).}. Let $\mu$ be an
$H$-invariant and ergodic measure. Then $\mu$ is homogeneous: i.e.
there is a closed group $L \leq G$ so that $\mu$ is an $L$-invariant
measure on a single $L$-orbit.
\end{Conjecture}

Given these conjectures, it is reasonable to expect that periodic
orbits of $H$ are nicely distributed in $\Gamma \backslash G$.
Indeed, as remarked by Margulis, the following follows from
Conjecture~\ref{bounded orbits conjecture} and the isolation results
of \cite{Cassels-Swinnerton-Dyer}:

\begin{Conjecture} \label{MargulisConjecture}
Let $n \geq 3$. For any compact $\Omega \subset \PGL_n(\Z) \backslash \PGL_n(\R)$
there are only finitely many periodic $H$-orbits contained in $\Omega$.
\end{Conjecture}

\subsection{Distribution and density of periodic orbits: statement of main results} \label{Resultstechnical}

\subsubsection{Periodic orbits and their approach to infinity, in the noncompact case.}

We obtain the following towards Conjecture \ref{MargulisConjecture}:
\begin{Theorem}\label{theorem about bounded orbits}
Let $n \geq 3$. For any compact $\Omega \subset \PGL_n(\Z) \backslash \PGL_n(\R)$ and any $\varepsilon > 0$ the total volume of all periodic $H$-orbits contained
in $\Omega$ of discriminant $\leq D$ is $\ll_{\varepsilon, \Omega} D ^ \varepsilon$.
\end{Theorem}

In contrast, for $n=2$ we have the following:

\begin{Theorem} \label{badcompactorbitsX2}
For every $\varepsilon>0$ there is a compact $\Omega \subset
\PGL_2(\Z) \backslash \PGL_2(\R)$ so that the total length of all
periodic $H$-orbits contained in $\Omega$ of discriminant $\leq D$
is $\gg_{\epsilon} D ^{1-\epsilon}$.
\end{Theorem}

Further,  even for $n \geq 3$, individual orbits may fail to become equidistributed, and, indeed,
spend a positive fraction of their time ``at $\infty$'', at least asymptotically --- see Section \ref{bad orbits section}.

\subsubsection{Density and distribution in the compact case.}
In the following results, let $\Gamma,G$ be as in (L-2) and $n \geq 3$.



Given a finite collection $Y = \{ x_iH: 1 \leq i \leq l \}$ of periodic orbits, we
define $$\disc(Y) := \max_{1 \leq i \leq l} \disc(x_iH), \ \ \vol(Y) :=
\sum_{i=1}^{l} \vol(x_{i}H).$$
Given such a collection $Y$, we let $\mu_{Y}$ be the sum of volume measures on $x_iH$, {\em normalized to have total mass $1$. }

\begin{Theorem}\label{theorem about density in the compact case}
Let $n \geq 3$ and let $\rho > 0$ be fixed. For each $j \geq 1$, let
$Y_j$ be a finite collection of periodic $H$-orbits satisying
$\disc(Y_j) \rightarrow \infty$ and $\vol(Y_j) \geq
\disc(Y_j)^{\rho}$. Suppose there is no periodic orbit of a group $H
< L < G$ (with both inclusions proper) containing infinitely many
periodic orbits belonging to $\bigcup Y_j$.
%
Then $\overline {\bigcup x_{j}H} = \Gamma \backslash G$.
\end{Theorem}

Our precise distribution result is stated in Theorem \ref{comprehensive equidistribution theorem}, but is rather complicated in general
because of the existence of intermediate subgroups. For now we simply state a corollary of that theorem:

\begin{Corollary}  \label{PrimeRankWellDistribution}(to Theorem \ref{comprehensive equidistribution theorem}).
Let $\Gamma,G$ be as in (L-2) with $n$ prime, and let $\rho > 0$ be fixed.

For each $j \geq 1$, let $Y_j$ be a finite collection of periodic $H$-orbits satisfying
$\disc(Y_j) \rightarrow \infty$ and
$\vol(Y_j) > \disc(Y_j)^{\rho}$.
Then any weak limit $\mu$ of the $\mu_{Y_j}$
satisfies
%
\begin{equation*}
\mu (B) \geq c_{n}.{\rho}.\mu_{\mathrm{Haar}} (B),
\end{equation*}
for any measurable set $B \subset X$; here $\mu_{\mathrm{Haar}}$ is
the $G$-invariant (Haar) probability measure on $\Gamma \backslash
G(\R)$ and $c_{n}$ is a positive constant which depends on $n$ only.
\end{Corollary}

\begin{Example}
A sequence $\{Y_{j}\}_{j\geq 1}$ to which this corollary applies is the following:
let $V_{1}<V_{2}<\dots<V_{j},\dots$ denote the ``volume spectrum'' of $X$  (i.e. the sequence of the volumes of all periodic $H$-orbits
  in increasing order), and let $Y_{j}$ be the collection of all orbits of
 volume $V_{j}$. It can be shown that the hypothesis of Corollary
   \ref{PrimeRankWellDistribution} is satisfied for any $\rho<1/2$ and $j$ large enough.
   \end{Example}

More generally, we can rephrase the corollary in the following way. Fixing a subset $U \subset \Gamma \backslash G(\R)$,
we call a periodic $H$-orbit $(U,\varepsilon)$-bad if the compact orbit spends
less than $\varepsilon \mu_{\mathrm{Haar}}(U)$ time inside $U$.
Then there is an constant $c_n > 0$, depending only on $n$,
so that the total volume of $(U,\varepsilon)$-bad orbits of discriminant $\leq D$
is $\ll_{\varepsilon, U}  D ^{\varepsilon c_n}$. In particular, the total volume of periodic orbits
that fail to intersect $U$ and have discriminant $\leq D$ is $\ll_{\varepsilon,U} D ^{\varepsilon}$.

\subsubsection{Linnik's principle.}

The key observation behind Theorem~\ref{theorem about bounded orbits} and Theorem~\ref{theorem about density in the compact case} is the following relation between total volume of a collection of periodic orbits in the relation to their discriminants and the entropy of any weak$^*$ limit.
In his book \cite{Linnik-book}, Yu. Linnik discusses several equidistribution problems, among them the equidistribution of packets of periodic geodesic trajectories on $\PGL_2(\Z) \backslash \PGL_2(\R)$.
Central to the proof of each of these equidistribution statements is a ``basic lemma'' which can be viewed as an implicit form of such a relation between volume and entropy,
and so we call this general phenomenon Linnik's Principle.

\begin{Theorem}  [Linnik's Principle] \label{thm:linnikprinciple} For $i \geq 1$,
let $Y_i$ be a collection of periodic $H$-orbits in $X = \Gamma \backslash G$
satisfying $\disc(Y_i) \rightarrow \infty$ and $\mathrm{vol}(Y_i) \geq \disc(Y_i)^{\rho}$, for some fixed $\rho > 0$.
Let $\mu _ i=\mu_{Y_{i}}$ be the probability measure associated with  $Y _ i$.  Suppose that
$\mu _ i \to \mu$ as $i \to \infty$ in the weak$^*$ topology
for some \emph{probability} measure $\mu$.
Then for any regular $h \in H$ there is an explicit $c_h>0$ (depending only on $h$) so that
\begin{equation*}
h _ \mu (h) \geq c_h \rho .\end{equation*}
\end{Theorem}

\subsection{Counterexamples to equidistribution} \label{bad orbits section}

Theorems \ref{theorem about bounded orbits} and \ref{theorem about density in the compact case} as well as Corollary~\ref{PrimeRankWellDistribution} fall far short of one possible candidate for an answer to the basic question of Section \ref{opening subsection}: that possibly periodic $H$-orbits of increasing volume which do not lie in any
periodic orbits of a bigger group\footnote {More formally, that there is no periodic orbits of a group $L$ with $H<L<G$ containing infinitely many of these periodic $H$-orbits.} become equidistributed in $\Gamma \backslash G$. This is not just because of technical difficulties; indeed, this plausible statement regarding equidistribution of $H$-orbits is \emph{false}, and the number of periodic $H$-orbits which fail to be equidistributed is asymptotically \emph{bigger than $D ^ \alpha$ for some strictly positive $\alpha$}:

\begin{Theorem} \label{badcompactorbits}
Let $n \geq 2$.
There is a sequence of periodic $H$-orbits $x_iH$ on the space $\PGL_n(\Z) \backslash \PGL_n(\R)$ and $\alpha, \delta >0$ such that
\begin{enumerate}
\item any weak$^*$ limit of the probability measures $\mu _ {x_iH}$ supported on these periodic $H$-orbits has total mass $<1-\delta$.
\item the total volume $\vol (\left\{ x _ iH : \disc (x_iH)< D \right\}) \gg D^ \alpha$.
\end{enumerate}
\end{Theorem}

\noindent This theorem is proved by an explicit construction (using a construction of Duke \cite{DukeCompositio}) suggested to us by P. Sarnak; a closely related construction in the special case of $n=2$ can be found in \cite{Sarnak-reciprocal-geodesics}.

One might hope escape of mass to infinity is the only possible obstacle to equidistribution of periodic $H$-orbits, and that things are nicer in the compact case. For example, one might hope that Corollary \ref{PrimeRankWellDistribution}
could be sharpened to give equidistribution, i.e. that the limiting measure $ \mu=\haarG$. Even this seems \emph{highly unlikely}.

While we do not have a counterexample in this particular setting (i.e. $G \cong \PGL (n, \R)$ and $\Gamma$ as in (L-2), with $n \geq 3$) we believe the following is likely: that for any periodic $H$-orbits $x_0H$ in $\Gamma \backslash G$,
there should exists collections $Y_i$ of periodic $H$-orbits with discriminants in the range $[\Delta_i,2 \Delta_i]$ and with $\vol (Y_i) \gg \Delta_i ^{\alpha}$ so that $\mu _ i$ converge weak$^*$ to a measure $\mu$ assigning positive measure to the given periodic $H$-orbit $x_0H$.
Moreover, we give an explicit construction of this type in the $S$-arithmetic context
in \cite{ELMV2}.

Besides showing the limit of what one may hope to prove regarding periodic torus orbits, these examples also point out to an important difference between torus actions and the action of groups generated by unipotents: \emph{linearization techniques \`a la Ratner, Dani-Margulis and Shah do not work for torus actions}, at least not on the level of individual periodic orbits, and moreover any substitute for these techniques will either have to assume much more or give much less\footnote {Indeed, to prove Theorems  \ref{theorem about bounded orbits} and \ref{theorem about density in the compact case} we use isolation results \cite{Cassels-Swinnerton-Dyer, Lindenstrauss-Barak} which can be viewed as a poor person's substitute for linearization techniques.}.

\subsection{A conjecture. Evidence from harmonic analysis.}
The construction of nonequidistributing periodic orbits described in Section \ref{bad orbits section} involve orbits $x_iH$ that have {\em very small volume relative to their discriminant}:
$\vol(x_iH)$ behaves as a polynomial in $\log \disc(x_iH)$.
We conjecture that this is the only source
of bad behavior:

\begin{Conjecture}\label{babyclass}
Fix $\rho > 0$. Let $G=\G(\R)$ be an $\R$-split real algebraic group, $\Gamma < G$ an arithmetic lattice, and $H$ a maximal $\R$-split torus. Let $x_i H$ be a sequence of periodic $H$-orbits satisfying $\vol(x_i H) \geq \disc(x_i H)^{\rho}$.
Then any weak limit of the measures $\mu_{x_i H}$ is algebraic.
\end{Conjecture}

We emphasize that the above Conjecture does not exclude the case of $\G =\PGL_2$.
At first sight, this may appear too optimistic, in view of the discussion of Section \ref{sec:higher rank rigidity}.
However, we expect the large torus orbits considered in Conjecture \ref{babyclass} to enjoy
a form of rigidity even in the case $\dim H=1$.
The reason for this is related to the concept of packets, mentioned in Section \ref{subsec:DSV}.
As remarked there, one may group the compact orbits into collections (packets) which admit
an action of an adelic group and have total volume $\disc^{1/2 +o(1)}$.

 It is likely that the large torus orbits considered in Conjecture \ref{babyclass}
will retain an action of a usable subgroup of this adelic group. This extra action provides a substitute
for a rank $2$ torus action.

Besides the results of \S \ref{Resultstechnical}, there exists
further evidence for Conjecture \ref{babyclass}, based on recent progress in the analytic theory of $L$-functions.
  Let us restrict to $\Gamma = \PGL_n(\Z), G = \PGL_n(\R)$.
In this context, the exponent $\alpha = 1/2$ is critical, because the total volume of all orbits
of discriminant $D$ is always $\gg_{\epsilon} D^{1/2-\epsilon}$.
In the case $n=2$, work of the third author with Harcos \cite{MichelSubconvex} implies
the conjecture for $\alpha$ slightly less than $1/2$. Moreover, the GRH
implies the conjecture for all $\alpha > 1/4$. These results
use (in addition to, respectively, the results of \cite{MichelSubconvex} and the GRH) the relationship between torus orbits and Rankin-Selberg $L$-functions,
established in the work of Waldspurger, Katok-Sarnak and Popa (see e.g.
\cite[Theorem 6.5.1]{Popa}).

It is worth emphasizing here that $\alpha = 1/4$ mark the limit of ``naive'' harmonic analysis.
One advantage of the ergodic theoretic methods used in this paper is that results such
as Corollary \ref{PrimeRankWellDistribution} yield nontrivial information even for very small $\alpha$.
In \cite{ELMV3}, we shall combine harmonic analysis and ergodic methods to obtain results for $n=3$.

\subsection{Analogy with $\times 2 \times 3$}

The dynamics of $x \mapsto 2x, x \mapsto 3x$ on $\R/\Z$ are, in many
ways, similar to the action of the diagonal group on
$\PGL_n(\R)/\PGL_n(\Z)$. Here the periodic points are exactly the
rationals $x =p/q$ with $(p,q) = 1=(q,6)$.

 Then the orbit $2^n . 3^m .x$ has two natural invariants: the size $g$ (i.e. the order of the group
$\langle 2,3 \rangle \subset (\Z/q\Z)^{\times})$ and the denominator $q$.  These
are analogues of the volume and discriminant of a compact $H$-orbit.   Again, the relation between
these two is a mystery, especially as far as how small $g$ can be; we do not even know that $g \geq 100 (\log q)^2$.

In their paper \cite{BGK} Bourgain, Glibichuk and Konyagin prove:
\begin{Theorem}[Bourgain-Glibichuk-Konyagin]
Let $\rho > 0$ be fixed.
Let $q$ be prime and let $G \subset (\Z/q\Z)^{\times}$
be a subgroup satisfying $|G| > q^{\rho}$. Then
for any $a \in (\Z/q\Z)^{\times}$, we have:
$$\left|\sum_{x \in G} e^{2 \pi i a x/q} \right| \leq |G| q^{-\delta}$$
where $\delta > 0$ depends only on $\rho$.
\end{Theorem}

In particular, if $q$ is a prime so that the order of $\langle 2, 3 \rangle \subset
(\Z/q\Z)^{\times}$ is large, then this theorem implies the periodic orbits
of $\times 2 \times 3$ with denominator $q$ are equidistributed individually,
as $q \rightarrow \infty$, in analogy with Conjecture \ref{babyclass}.
The restriction that $q$ be prime has been lifted in
\cite{ChangBourgain} and subsequent work.

The analogy with $\times 2 \times 3$ is also useful to illustrate, in an explicit context, many of the
ideas that are used in this paper. In this vein, we present the analogue of Corollary \ref{PrimeRankWellDistribution}
as well as a sketch of its proof.

\begin{Proposition} Fix $\rho >\rho'> 0$ and a subinterval $J \subset \R/\Z$.
 Let $S \subset q^{-1} \Z/\Z$
be invariant under $x \mapsto 2x, x \mapsto 3x$ and so that $|S| > q^{\rho}$.
Then, for sufficiently large $q$, we have
$$ \frac{|S \cap J|}{|S|} \geq \rho' \mathrm{length}(J)$$
\end{Proposition}

  Equivalently: let $S_i \subset q_i^{-1} \Z/\Z$ with $|S_i| > q_i^{\rho}$
be invariant under $x \mapsto 2x, x \mapsto 3x$.
Let $\mu_i$ be the corresponding normalized measure\footnote{That is, $\mu_i(f) = \frac{1}{|S_i|}
\sum_{x \in S_i} f(x)$ for a continuous function $f$ on $\R/\Z$.}. Then any weak limit of the $\mu_i$
dominates $\rho. \mu_{\mathrm{Lebesgue}}$.

\proof
We make free use of the notion of entropy.\footnote{See \eqref{entropydef} for the definition and
\cite{Katok-book} for background. We will use the fact -- subadditivity of entropy -- that $H_{\mu}(\mathcal{P} \vee \mathcal{Q})
\leq H_{\mu}(\mathcal{P}) + H_{\mu})(\mathcal{Q})$ for two partitions $\mathcal{P}, \mathcal{Q}$.}
Let $\mathcal{P}$ be the partition of $\R/\Z$ into $[0,1/2) \cup [1/2,1)$;
let $\mathcal{P}^{(n)} := \mathcal{P} \vee [2]^{-1} \mathcal{P} \vee \dots [2^{n-1}]^{-1} \mathcal{P}$.
Here $[m]^{-1} \mathcal{P}$ denotes the partition into $\{x: mx \in [0,1/2)\}$ and  its complement.

Let $n_i$ be minimal so that $2^{n_i} > q_i$. Then any two elements of $S_i$ lie in distinct parts
of the partition $\mathcal{P}^{(n_i)}$. From this it follows that $H_{\mu_i}(\mathcal{P}^{(n_i)})
> \log(q_i)$. By the subadditivity of entropy,
$$ H_{\mu_i}(\mathcal{P}^{(n_i)}) \leq  H_{\mu_i}(\mathcal{P}) + H_{\mu_i}([2]^{-1} \mathcal{P}) + \dots
+ H_{\mu_i}([2^{n_i-1}]^{-1} \mathcal{P})$$
Because $\mu_i$ is invariant under $x \mapsto 2x$, the summands on the right-hand side are equal.
Therefore, $H_{\mu_i}(\mathcal{P}) \geq \frac{\log |S_i|}{n_i}$. Passing to the limit, we conclude
that any weak limit $\mu$ of the $\mu_i$ satisfies:
\begin{equation} \label{basicentropy} H_{\mu}(\mathcal{P}) \geq \rho \log 2. \end{equation}

By a simple variation of this argument, one verifies that
$H_{\mu}(\mathcal{P}^{(k)}) \geq k \rho \log 2$. This implies that
the entropy of the transformation $x \mapsto 2x$ w.r.t. $\mu$ is
$\geq \rho \log 2$.

We may now apply the following theorem of Rudolph \cite{Rudolph-2-and-3}:
A probability measure on $\R/\Z$, ergodic and invariant for $\times 2 \times 3$,
which has positive entropy w.r.t. $\times 2$, is Lebesgue measure.

To apply it, we decompose $\mu$ into ergodic components with respect to $\times 2 \times 3$.
Because any probability measure on $\R/\Z$ has entropy $\leq \log 2$ w.r.t. $\times 2$,
it follows that $\mu$ must dominate a measure $\mu'$ with total mass $\geq \rho$
and so that almost every ergodic component of $\mu'$ has positive entropy w.r.t. $\times 2$.
By Rudolph's theorem, $\mu'$ dominates $\rho . \mu_{\mathrm{Lebesgue}}$.
\qed

The key points in the above proof were: the fact that elements of $S_i$
were well-separated at a small scale (namely $q_i^{-1}$); the use of subadditivity of entropy;
and Rudolph's theorem. The proof of most of the results of this paper, including Corollary \ref{PrimeRankWellDistribution}, will
use similar ideas in the $\Gamma \backslash G$ context.
The well-separated property is established in Proposition \ref{baby}, and
we replace the use of Rudolph's theorem
with the results of \cite{Einsiedler-Katok-Lindenstrauss}.

\begin{Acknowledgements} The present paper is part of a project that began on the occasion of the AIM workshop ``Emerging applications of measure rigidity" on June 2004 in Palo Alto.
 It is a pleasure to thank the American Institute of Mathematics,
 as well as the organizers of the workshop.
 We would also like to thank Peter Sarnak for his suggestions, comments, and for his constant encouragement
 during the elaboration of this paper.

  While working on this project the authors visited the following institutions: Princeton University (Ph. M.), the Center of Mathematical Sciences, Zhejiang University  (Ph. M. and A. V.),
  the Institut des Hautes Etudes Scientifiques  (Ph. M. and A. V.). The support of these
institutions are gratefully acknowledged.
\end{Acknowledgements}
\section{The description of periodic orbits} \label{sec:gn}
Let notation be as in Section \ref{subsec:DSV}.

Our aim in this section is to discuss the parameterization of periodic Cartan orbits and to attach
to each orbit an integral invariant, the {\em discriminant}.  The point is that periodic orbits
will be parameterized by rational points on a certain variety (namely, the variety
of tori inside $\G$) and the discriminant will be the denominator of the associated point.

\begin{enumerate}
\item Two periodic orbits of discriminant $\leq D$ cannot be too close:
if $\Omega \subset  G$ is compact,
and $g_i \in \Omega$ are such that $g_1 H \neq g_2 H$, and $\Gamma g_i H$ are periodic,
then the distance between $g_1, g_2$ is $\gg_{\Omega} D ^{-1}$. (See Proposition \ref{proposition discrete}, which shows a sharper form of this statement).
\item A periodic orbit of discriminant $\leq D$ cannot go too close to $\infty$ -- this is of relevance
only if $\Gamma \backslash G$ is noncompact (Proposition \ref{proposition bounding periodic orbits}.)
\item The volume of a periodic orbit of discriminant $D$ is bounded, up to constants,
between $\log D$ and $D ^{c}$, for $c >0$ (Proposition \ref{proposition about regulator}).
\item  The number of periodic Cartan orbits of discriminant $\leq D$ is bounded, up to constants, between $D ^{A}$ and $D ^{B}$, for $A, B >0 $ (Proposition \ref{proposition about number of orbits}).
\end{enumerate}

A serious omission in the present account will be that we will avoid any explicit discussion of packets,
the natural equivalence relation on periodic orbits. This will be discussed in the paper \cite{ELMV3}; see also Remark \ref{Packet} in the present paper.

The paramaterization of periodic orbits that we present here is essentially well-known; see e.g. \cite{Oh-finiteness}. The primary goal of this section is, rather, to define the discriminant and show that it crudely controls
the dynamics of the orbit, as discussed above.

In Section \ref{Section algebras} we will specialize some of the discussion presented here to
algebraic groups defined by central simple algebras and to the case of $\PGL_n$.

\subsection{Notation}
As a general rule, we will use {b}oldface letters for algebraic varieties, algebraic groups and so forth.

Let $\G$ be a semisimple group over $\Q$ and put $G=\G(\R)$. We will assume throughout this document
that $\G$ is split over $\R$. While the techniques and ideas extend to the general case,
this assumption allows for the most elegant and coherent treatment.

Let
$\Gamma$ be an arithmetic subgroup of $G$ that is commensurable with $\G (\Q) \cap \GL (n, \Z)$ for some fixed embedding $\rho$ of $\G$ into $\GL (n)$ defined over $\Q$.  Let $X = \Gamma
\backslash G$ be the associated homogeneous space. Let $\mathfrak{g} =
\Lie(\G)$ be the Lie algebra of $\G$; thus $\mathfrak{g}$ is a Lie
algebra over $\Q$. We fix once and for all a lattice
$\mathfrak{g}_{\Z}$ in $\mathfrak{g}$ and a $\G$-invariant nondegenerate bilinear
symmetric form $B(\cdot, \cdot)$ defined by
$$B(X,Y) = \mathrm{Tr}(\rho(X) \rho(Y))$$
for some fixed faithful representation $\rho$.
We assume that these
choices have the properties that: $\mathfrak{g}_{\Z}$ is stable by
the adjoint action of $\Gamma$, $B(\mathfrak{g}_{\Z},
\mathfrak{g}_{\Z} ) \subset \Z$, and $[\mathfrak{g}_{\Z}, \mathfrak{g}_{\Z}] \subset \mathfrak{g}_{\Z}$. We can always find such a $\mathfrak{g}_{\Z}$ (once one has a $\Gamma$-stable lattice,
a suitable integral multiple of it will satisfy the latter two properties).

Let $H$ be an $\R$-split Cartan subgroup of $G$ (i.e.\ $H$
is the centralizer of a split Cartan subalgebra of $\g \otimes \R$.)
The chosen data fix a Haar measure on $H$: the form $B$ restricted to the Lie algebra $\h = \mathrm{Lie}(H)$ is nondegenerate and we shall use the Haar measure on $H$ determined by a top differential form on $\h$ self-dual with respect to $B$.

Fix an Euclidean norm $\norm {\cdot}$ on $\mathfrak g$; this also induces an Euclidean norm on $ \wedge ^ \ell \mathfrak g$ for any $\ell$.
The choice of Euclidean norm on $\mathfrak g$ determines a left invariant Riemannian metric $d(\cdot, \cdot)$  on $G$ which induces also a metric on $X = \Gamma \backslash G$.

Recall (see, e.g. \cite[Remark 3.6]{Tomanov-Weiss}) that for any $R > -1$ the set
\begin{equation}\label{omegadef}
\Omega(R) = \left\{ \Gamma g: \norm {\operatorname{Ad} (g ^{-1}) v} \geq R ^{-1} \text{ for all nonzero $ v \in \mathfrak g _ \Z$} \right\}
\end{equation}
is compact with $\bigcup_ R \Omega (R) = \Gamma \backslash G$. It
can easily be shown, e.g. using reduction theory \footnote{Indeed,
fix a maximally $\Q$-split torus $\mathbf{T}_0 \subset \mathbf{G}$,
and let $A_0 = \mathbf{T}_0(\R)$. Take $x \in \Omega(R)$. It has by
reduction theory a representative $\delta. a. \kappa \in G$, where
$a \in A_0$, $\delta$ belongs to a finite subset of $\G(\Q)$, and
$\kappa$ belongs to a fixed compact subset of $G$. But then $\inf_{v
\in \mathfrak{g}_{\Z}} \|\Ad(g)^{-1} v \| $ is bounded above and
below by multiples of $\|\Ad(a ^{-1})\|$, because the eigenspaces of
$\mathrm{Ad}(a ^{-1})$ are $\Q$-rational subspaces of
$\mathfrak{g}$.  Moreover, $\|\Ad(a ^{-1})\|$ is bounded above and
below by multiples of $\|\Ad(g ^{-1})\|$. Finally, note that
$\|\Ad(g)\|$ is bounded from above and below by multiples of
$\sup_{\|u\|,\|v\|=1}|B(\Ad(g)u,v)|$ and so also by multiples of
$\|\Ad(g^{-1}\|$.} that
every $x \in \Omega (R)$ can be represented as $x = \Gamma g$ for some $g$ with $\norm {\operatorname{Ad} (g)} \ll R$; here $\| \operatorname{Ad}(g)\|$ denotes the operator norm
of $\operatorname{Ad}(g)$.
with respect to the fixed Euclidean norm on $\mathfrak{g}$.

Finally, in what follows, we will allow the implicit constant in constructions such as $O(\cdot), \ll, \gg$
to depend on the data fixed above (in particular $\G, \Gamma, \mathfrak{g}_{\Z}$).

\subsection{The parameterization of periodic orbits.} \label{subsec:gen}
We discuss in this section how to parameterize periodic orbits of $H$ on $X$.

Given a periodic orbit $\Gamma g H$ of $H$ on $X$, the stabilizer
in $H$ of any point in this orbit, i.e. $H \cap g ^{-1} \Gamma g$, is necessarily
a cocompact lattice in $H$. The Zariski closure of $\Gamma \cap g H g ^{-1}$
is then a torus $\T$ defined over $\Q$ with the property that $\T(\R)=g  H g ^{-1}$.
The fact that $\Gamma \cap \T(\R)$ is an arithmetic cocompact subgroup of $\T(\R)$ assures us that
$\T$ is $\Q$-anisotropic. \footnote{A character $\chi$ defined on $\T$ over $\Q$ would map the integer points $\T (\Z)$ to a subgroup $\chi (\T (\Z))$ of $\Q ^ \times$ whose elements allow a common denominator, clearly this subgroup must be trivial and so the same holds for $\chi$ since $\T (\Z)$ is a cocompact subgroup of $\T(\R)$.}



\begin{Proposition}[Basic correspondence]  \label{Basic bijection}
There is a canonical bijection between
\begin{enumerate}
\item [(1)] periodic $H$-orbits $\Gamma g H$ on $\Gamma \backslash G$, and
\item [(2)]
$\Gamma$-orbits on pairs $(\T, g H)$.  Here $\T$ is a maximal $\Q$-torus
that is anisotropic and $\R$-split, and $g H \in G/H$ is such that $g H g ^{-1} = \T(\R)$, and $\gamma \in \Gamma$
acts via $(\T, gH) \mapsto (\gamma \T \gamma ^{-1} , \gamma g H)$.
\end{enumerate}
\end{Proposition}


\begin{proof}
The correspondence has already been indicated: to a periodic $H$-orbit $\Gamma g H$
we associate $(\T, gH)$, where $\T$ is the Zariski closure of $\Gamma \cap g H g ^{-1}$.
In the reverse direction, we associate to $(\T, g H)$ the orbit $\Gamma g H$;
it is compact because $\Gamma \cap \T(\R)$ is a cocompact lattice in $\T(\R)$.
\end{proof}

We remark that $\T$, which is necessarily $\R$-split, determines the periodic $H$-orbit up to finitely many possibilities.
Indeed, each $ \mathbf T$ corresponds to precisely one periodic $N _ G (H)$-orbit, and so to at most $\absolute {H \backslash N _ G (H)}$ many periodic $H$-orbits.

\subsection{The discriminant of a periodic orbit.} \label{discriminantsection}
We would like to attach to each periodic orbit $x_0 H$ a positive integral invariant,
the {\em discriminant}.

We will do this by means of the correspondence of Proposition \ref{Basic bijection}: we will assign a positive
integer to each maximal $\Q$-torus contained in $\G$, in a $\Gamma$-invariant fashion.
This will be done by interpreting $\Q$-torii as
rational points on some variety, and interpreting the discriminant as the ``least common denominator''
in the coordinates of these points.\footnote{We remark that this ``discriminant'' is not a canonical construction, e.g.\ it depends on the choices of the lattice $\g_\Z$ and the bilinear symmetric form $B (\cdot, \cdot)$.  However, as we will see later we have set it up so that it agrees with the notion of discriminant for a number field or order (assuming the right choices have been made).}

Let $r$ be the rank of $\mathbf{G}$.
The bilinear form $B$ on $\mathfrak g$ gives rise to a bilinear form $B _ {\wedge}$ on $\wedge ^ r \mathfrak g$ determined by
\begin{equation*}
B _ \wedge (x _ 1 \wedge \dots \wedge x _ r, y _ 1 \wedge \dots \wedge y _ r) = \det (B (x _ i, y _ j)) _ {i, j = 1} ^ r
.\end{equation*}
Set $V = \wedge ^ r \mathfrak g \otimes \wedge ^ r \mathfrak g$, $V_\Z = \wedge ^ r \mathfrak g _ \Z \otimes \wedge ^ r \mathfrak g _ \Z$.

For $g \in G$, we will often use the notation $\Ad(g)$ to denote the natural action
of $g \in G $ on any of $\mathfrak{g}$ (the adjoint action), $\wedge ^ r \mathfrak{g}$
and $V$.

For any maximal $\R$-split Cartan subalgebra $\mathfrak t < \mathfrak g$ we define $\iota (\mathfrak t) \in V$ by taking any $w \in \wedge ^ r \mathfrak t$ and setting
\begin{equation} \label{iotadef}
\iota (\mathfrak t) = \frac {w \otimes w }{ B _ \wedge (w, w)}
.\end{equation}

Since the form $B$ restricted to an $\R$-split Cartan subalgebra
$\mathfrak t$ is positive definite, for any nonzero $w \in \wedge ^
r \mathfrak t$ we have that $B _ {\wedge ^ r \mathfrak g}(w,w) > 0$.
If $\mathfrak t$ is the Lie algebra of a $\Q$-torus, $\mathfrak t _
\Z = \mathfrak t \cap \mathfrak g _ \Z$ is a lattice in $\mathfrak
t$. Choose a basis $e _ 1, \dots, e _ r$ of $\mathfrak t _ \Z$; then
$w = e _ 1 \wedge \dots \wedge e _ r$ is a primitive element of
$\wedge ^ r \mathfrak g _ \Z$, $B _ \wedge (w, w)$ is a positive
integer, and is the smallest integer $n$ such that $n \iota
(\mathfrak t) \in V _ \Z$.

The map $g \mapsto g \mathfrak h g ^{-1}$ gives a bijection between
$G / N _ G (H)$ and $\R$-split Cartan subalgebras of $\mathfrak g$,
hence we can view $\iota$ as a map from $G / N _ G (H)$ to $V$. If
$\Gamma g H$ is a periodic $H$ orbit $g \mathfrak h g ^{-1}$ is the
Lie algebra of an $\R$-split $\Q$-torus of $G$.

\begin{Definition}
Let $\mathfrak t$ be the Lie algebra of a maximal $\R$-split $\Q$-torus of $G$;
We define the {\em discriminant} $\disc(\mathfrak t) = \min \{ n \in \N: n  \iota (\mathfrak t) \in V_{\Z}\}$.
We define
the {\em discriminant} of the periodic orbit $\Gamma g H$ to be the discriminant
of $g \mathfrak h g ^{-1}$.
\end{Definition}

We remark that $\disc(\Ad (\gamma) \mathfrak t) = \disc(\mathfrak
t)$ for $\gamma \in \Gamma$, by virtue of the assumption that
$\g_{\Z}$ is invariant under the adjoint action of $\Gamma$, and
hence the discriminant of a periodic orbit is well-defined.

\subsection {Discriminant, discreteness and related properties of closed orbits}

The definition of the discriminant of periodic $H$-orbits, and its interpretation in terms of integral points on hypersurfaces, imply several useful facts. The fact that integral points are spaced by $\geq 1$
translates to the following  basic proposition, showing that
two periodic $H$-orbits of small discriminant cannot be too close. Though simple, this observation is absolutely crucial to our approach.

\begin{Proposition} \label{baby} (Discreteness of periodic orbits -- weaker form.)
Let $\Omega \subset G$ be a fixed compact subset of $G$.
Suppose $g_1, g_2 \in \Omega$ are so that $\Gamma g_i H$ are periodic orbits of discriminant $\leq D$ with $g_1 N_G(H) \ne g_2 N_G(H)$.
Then $d(g_1, g_2) \gg_{\Omega} D ^{-2}$.
\end{Proposition}
Note that we do not assume that $\Gamma g_1 N_G(H) \ne \Gamma g_2 N_G(H)$, i.e.\ the proposition also describes how close a particular periodic orbit can come to itself.
\begin{proof}
Let $\mathfrak{t}_i = \Ad(g_i) \mathfrak{h}$. Then
$\| \iota(\mathfrak{t}_1) - \iota(\mathfrak{t}_2)\| \gg D ^{-2}$, as is clear from \eqref{iotadef}.
But the map $g \mapsto \iota(\Ad(g) \mathfrak{h})$ is a smooth map from $\Omega$ to $\wedge ^{r} \mathfrak{g}$; in particular, it cannot increase distances by more than a constant factor depending on $\Omega$.
\end{proof}

Actually we prefer to give a slightly sharper version of this result, which both explicates
the dependence on $\Omega$ and gives better dependence on $D$.

\begin{Proposition} \label{proposition discrete} (Discreteness of periodic orbits -- sharper form.)
Suppose $g_1, g_2$ satisfy $\|\Ad(g_i)\| \leq R$ and are such that  $\Gamma g_i H$ are periodic orbits of discriminant $D_i$ respectively  with $g_1 N_G(H) \ne g_2 N_G(H)$
Then
\begin{equation*}
d(g_1, g_2) \gg \begin{cases} R ^{-r} D ^ {-1/2}& D=D_1=D_2\\
R ^{-2r} (D _ 1 D _ 2) ^{-1/2}&D_1 \neq D_2
\end{cases}
.\end{equation*}
\end{Proposition}

%

\begin{proof}
Let $\Gamma g_i H$ be periodic, $\mathfrak t _ i = g _ i \mathfrak h g _ i ^{-1}$ and $D_i = \disc (\Gamma g _ i H)$.
Let $w_0 \in \wedge ^ r \mathfrak h$ be such that $B _ \wedge (w _ 0, w _ 0) = 1$ (this fixes $w _ 0$ up to sign), and set $w _ i = \operatorname{Ad} (g _ i ) [w _ 0]$. Since $g _ 1 N _ G (H) \neq g _ 2 N _ G (H)$ we have that $w _ 1 \neq w _ 2$.

It follows from the definition of discriminant that $D _ i ^{1/2} w _ i \in \wedge ^ r \mathfrak g _ \Z$.
Indeed, if $e _ 1, \dots, e _ r$ are a basis of $\mathfrak t _ i \cap \mathfrak g _ \Z$ then $w= e _ 1 \wedge \dots \wedge e _ r \in \wedge ^ r \mathfrak g _ \Z$ is an element of $\wedge ^ r \mathfrak t$ with $B_\wedge (w,w)=D_i$. On the other hand $D _ i ^{1/2} w _ i$ has the same property; hence
$D _ i ^{1/2} w _ i = \pm w \in \wedge ^ r \mathfrak g _ \Z$.

Consider first the case $D _ 1 = D _ 2 = D$.  Then, since $D ^{1/2} w _ 1$ and $D ^{1/2} w _ 2$ are two distinct ``integral points'', we have $\norm {w _ 1 - w _ 2} \gg D ^ {- {\tfrac{1}{2}}} $, i.e.
$\| \Ad(g_1) w_0 - \Ad(g_2) w_0 \| \gg D ^{-1/2}$.

The operator norm of $\Ad(g_1)$ acting on $\mathfrak{g}$ is $\leq
R$; therefore, the operator norm of $\Ad(g_1)$ acting on $\wedge
^{r} \mathfrak{g}$ is $\leq R^r$, and we conclude that
\begin{equation*}
\| w_0 - \Ad(g_1 ^{-1} g_2) w_0 \|
\gg
R ^{-r} D ^ {- 1/2}
.\end{equation*}
This shows that $d(g_1,g_2) = d(e, g_1 ^{-1} g_2) \gg R ^{-r} D ^ {- 1/2}$.

Now consider the case $D _ 1 \neq D _ 2$.
Then $D _ 1 ^{1/2} w _ 1$ and $D _ 2 ^{1/2} w _ 2$ are two nonproportional integer points, hence they span a parallelogram of area $\gg 1$ and so
\begin{equation} \label{parallelogram equation}
\norm {D _ 2 ^ {1/2} (w _ 2 - w _ 1)} \norm {D _ 1 ^ {1/2} w_1} \gg 1
.\end{equation}

As before, the operator norm of $\Ad(g_1)$ on $\wedge ^ r
\mathfrak{g}$ is $\leq R^r$, so $\norm {w_1} = \norm
{\operatorname{Ad} (g_1) w_0} \ll R ^{r}$ and so from
\equ{parallelogram equation} we get that
\begin{equation*}
\norm {w _ 1 - w _ 2} \gg R ^ {- r} (D _ 1 D _ 2) ^ {- 1 / 2}
.\end{equation*}
Proceeding again as in the previous case, we see that
$d(g _ 1, g _ 2) \gg R ^{-2r} (D_1D_2) ^ {- 1/2}$.
\end{proof}

The discriminant also bounds how high a periodic $H$ orbit can penetrate the cusps of $\Gamma \backslash G$. Recall \eqref{omegadef} for the definition of $\Omega(\cdot)$.


\begin{Proposition}\label{proposition bounding periodic orbits}
Let $\Gamma g H$ be a periodic $H$-orbit of discriminant $D$. Then
$\Gamma g \in \Omega (c D ^{\dim \mathfrak{g}/2})$ for some
constants $c> 0$ (independent of $g$).
\end{Proposition}

\begin{proof}
The idea is that $\mathfrak{h} \cap \mathrm{Ad}(g ^{-1}) \mathfrak{g}_{\Z}$
has (by definition of discriminant, and lattice reduction) a basis consisting of vectors that are not too long. From this, we need to finesse that
the same is actually true for $\mathrm{Ad}(g ^{-1}) \mathfrak{g}_{\Z}$. But if it contained short vectors,
they would generate a nilpotent subalgebra $\mathfrak{n}$ normalized by $\mathfrak{h}$;
this would mean that $\mathrm{Ad}(g) \mathfrak{n}$ would be a $\Q$-rational
nilpotent algebra normalized by $\mathfrak{t} = \mathrm{Ad}(g) \mathfrak{h}$, a contradiction,
for the associated torus is anisotropic.

Let $e_1, \dots, e_r$ be a basis for $\mathfrak{t} \cap \g_{\Z}$.
Since $B (e _ i, e _ j) \in \Z$ for every $i,j$ and $\det (B (e _ i, e _ j )) = D$, lattice reduction shows that we can choose the basis $e _ i$ so that $B(e _ i, e _ i) \leq c_0 D$, the constant $c_0$ depending only on $r$. Since $e_i' := \operatorname{Ad} (g) e_i \in \mathfrak h$, and on $\mathfrak h$ the form $B(\cdot, \cdot)$ is positive definite, this implies that
\begin{equation} \label{equation about conjugating by g}
\norm {e_i'}   \leq c_1' B(e_i', e_i') ^{1/2}  \leq c_1 D ^{1/2} \qquad \text{for $1 \leq i \leq r$}
\nonumber
.\end{equation}
Therefore, $\mathfrak{h} \cap \mathrm{Ad}(g ^{-1}) \mathfrak{g}_{\Z}$ has a basis
of vectors of length $\leq c_1 D ^{1/2}$. We may assume $c_1 > 1$, increasing it if necessary.

On the other hand,
there is some $c_2 \geq 1$ so that
\begin{equation*}
\norm {[v,w]} \leq c _ 2 \norm {v} \norm {w}
.\end{equation*}
For $i \geq 0$, let $\mathfrak{g}_i$ be the subspace spanned by $X \in \mathrm{Ad}(g ^{-1}) \mathfrak{g}_{\Z}$
with $\| X \| \leq  c_2 ^{-1} (c_2 c_1 D ^{1/2})^{-i}$. Then for $i \geq 1$ we have
$[\mathfrak{h}, \mathfrak{g}_i] \subset \mathfrak{g}_{i-1}$;
moreover, for any $i, j \geq 1$ we have $[\mathfrak{g}_i, \mathfrak{g}_j] \subset \mathfrak{g}_{i + j }$.
So there is
$i \leq \dim(\mathfrak{g})$ so that $\g_i = \g_{i+1}$.
Then $\mathfrak{n} := \mathrm{Ad}(g) \g_i$ is, {\em if nonempty}, a $\t = \Ad(g) \mathfrak{h}$-stable nilpotent subalgebra of $\mathfrak{g}$ defined over $\Q$.

%


It follows that $\mathfrak{t}$ normalizes the nilpotent Lie algebra
$\mathrm{Ad}(g) \mathfrak{n}$, which is defined over $\Q$.  But this
is a contradiction unless $\mathfrak{n}$ is trivial, for $\T$ is
anisotropic.\footnote{Over an algebraic closure, $\t \oplus
\mathfrak{n}$ is solvable and so contained in a Borel subalgebra,
and so all the roots of $\T$ on $\mathfrak{n}$ are positive in a
suitable system. This means that the determinant of the action of
$\T$ on $\mathfrak{n}$ cannot be trivial, and so defines a
nontrivial character of $\T$.}

We conclude that $\Ad(g ^{-1}) \mathfrak{g}_{\Z}$ had no nonzero elements of norm $\leq c_3 D ^{-\dim \mathfrak{g}/2}$.
\end{proof}

The exponent here could be considerably improved. It is easy to see, for instance,
that the exponent $\dim \mathfrak{g}/2$ could be replaced by $1/2$ for $\G = \PGL_2$ (see also
Corollary \ref{Cor:Explication}).  This is indeed sharp:

\begin{Example}
Let $G = \PGL_2(\R)$, $\Gamma = \PGL_2(\Z)$, $H$ the group of diagonal matrices and $d$ a square free positive integer. Then for $g = \begin{pmatrix} 1 & 1\\
\sqrt d& - \sqrt d \end{pmatrix}$ the point $\Gamma g$ is periodic under $H$, of discriminant a constant multiple of $d$ (assuming one makes the obvious choices for $\mathfrak g _ \Z$ etc.\footnote{These (unilluminating) choices are set up for $\PGL_n(\Z) \backslash \PGL_n(\R)$ in Section \ref{dynamics section}.}) Take $v = \begin{pmatrix} 0& 1\\
0& 0 \end{pmatrix} \in \mathfrak g _ \Z$. Then
\begin{equation*}
\operatorname{Ad} (g) v = \frac {1 }{ 2}\begin{pmatrix} 1/\sqrt d& 1 / \sqrt d\\
- 1 / \sqrt d& - 1 / \sqrt d
\end{pmatrix}
.\end{equation*}
This implies $g \in \Omega(c_0 d ^{1/2})$ for a suitable absolute constant $c_0$.
\end{Example}

We recall the following useful fact about $H$-orbits:

\begin{Theorem} [Tomanov and Weiss {\cite [Thm 1.3]{ Tomanov-Weiss}}]
\label{theorem of Tomanov and Weiss}
There is an $R_0$ so that for every $x \in \Gamma \backslash G$ we have that $\Omega (R_0) \cap xH \neq \emptyset$.
\end{Theorem}

There is a relation between the discriminant and the volume of a periodic $H$-orbit. This connection is explicated in the context of central simple algebras in \S\ref{Section algebras}. Here we give the following general (and less accurate) bounds:

\begin{Proposition} \label{proposition about regulator}
Let $ \Gamma g_0 H$ be periodic $H$-orbit. Then
\begin{equation} \label{regulator and discriminant inequality}
\log \disc (\Gamma g _ 0 H) \ll \vol (\Gamma g _ 0 H)
\ll \disc (\Gamma g _ 0 H) ^ c
.\end{equation}
for some $c$.
\end{Proposition}

\begin{proof}
Let $x_0 = \Gamma g_0$ be a point with $x_0 H$ a periodic $H$-orbit.
By Theorem~\ref{theorem of Tomanov and Weiss}, we may assume $g_0$ is in some fixed compact subset $\Omega \subset G$. Since $g_0 \in \Omega$, there are generators $h _ 1, \dots, h _ r$ of $g_0 ^{-1} \Gamma g_0 \cap H$ (possibly up to finite index) so that
\begin{equation*}
\vol (x_0 H) = \vol (g_0 ^{-1} \Gamma g_0 \backslash H) \gg \max_ i d(e, h_i)
.\end{equation*}
It follows that, if one chooses a basis of $\mathfrak{g}$ belonging
to the lattice $\mathfrak{g}_{\Z}$, then $\Ad(\gamma _ i) = \Ad(g_0 h _ ig_0 ^{-1})$
is represented (w.r.t. this basis) by an integral matrix all of whose coordinates are $\exp (O (\vol (x_0 H)))$.

The Lie algebra $g_0 \mathfrak h g_0 ^{-1}$ is precisely the
subspace fixed by the action of all $\gamma _ i$. Thus there
is\footnote{Given integer matrices in $M_n(\Z)$, all of whom have
all matrix entries $\leq N$, it is easy to see -- Siegel's lemma --
that their common kernel, if nonempty, contains an element of length
$\ll N ^{c(n)}$, where $c(n)$ depends only on $n$.} a nonzero $w \in
\wedge ^ r \mathfrak{g}_{\Z} \cap  \wedge ^ r g _ 0 \mathfrak h g _
0 ^{-1}$ with $\norm {w} \ll \exp (O (\vol (x _ 0 H)))$. One
concludes:
\begin{equation*}
\log \disc (x_0 H) \leq \log B _ {\wedge ^ r \mathfrak g} (w,w) \ll \vol (x_0 H)
.\end{equation*}
This proves the lower bound in \equ{regulator and discriminant inequality}.

To prove the upper bound in \equ{regulator and discriminant
inequality}, let $\mathfrak S \subset G$ be such that $G = \Gamma
\mathfrak S$ and $\mathfrak S$ has finite Haar measure. An example
of such a set is furnished by Siegel domains (see \cite [Sec. 4]{
Borel-Harish-Chandra}). Let $\mathfrak S (R) = \{g\in \mathfrak S:
\Gamma g\in \Omega (R)\}$. These Siegel domains have the following
additional property which will be useful for us (and is easily
verified directly from their definition): for a suitable $C,
\beta>0$
\begin{equation*}
\haarG (B_{R ^{-\beta}}(\mathfrak S (R))) < C \qquad \text{for all $R >0$}
\end{equation*}
where $B_{\delta}(S) = \left\{ g \in G: d(g,S)<\delta \right\}$.

Let $D = \disc (x_0 H)$. By Proposition~\ref{proposition bounding
periodic orbits}, the orbit $x _ 0 H \subset \Gamma\mathfrak S (cD ^
{\dim \g/2})$, and so we can write $ x_0 H \subset G$ as a (finite
or countable) disjoint union $\bigcup_ {i \in J} \Gamma S_i$ with
each $S_i \subset g_i H \cap \mathfrak S (c D ^ {\dim \g/2})$ for
some $g_i \in G$ lying in distinct cosets of $H$. Let $J' \subset J$
be such that $g_i N_G(H) \neq g _ j N_G(H)$ for every $i,j \in J'$
and such that
\begin{equation*}
\sum_ {i \in J '} \vol (S_i)
\end{equation*}
is maximal among all possible choices of $J '$. Then $\sum_ {i \in J '} \vol (S_i) \nope \gg \vol (x_0 H)$.

By Proposition~\ref{proposition discrete}, for every pair $i,j \in J '$ with $i \ne j$ we have that
\begin{equation*}
B _ \delta (S _ i) \cap B _ \delta (S _ j) = \emptyset \qquad \text
{for $\delta \leq 2 c_4 D ^ {-c_5}$}
\end{equation*}
$0<c_4<1$ and $c_5\geq  \beta$ being some constants independent of
$x_0$. Then
\begin{align*}
1 & \gg \haarG (B_{c_4D ^{-c_5}}(\mathfrak S (cD ^{\dim \g/2}))) \\
&\geq \sum_ {i \in J '} \haarG (B_{c_4 D ^{-c_5}} (S_i))\\
& \gg \sum_ {i \in J '} D ^{-c'} \vol (S_i) \gg D ^{-c'} \vol (x_0 H)
.\end{align*}
\end{proof}

\begin{Proposition}\label{proposition about number of orbits}
Let $N(D)$ be the number of periodic orbits with discriminant less than $D$.
Then there exists $c_1,c_2$ such that $D ^{c _ 1} \ll N(D) \ll D ^ {c _ 2}$.
\end{Proposition}

\begin{proof}
[First proof of upper bound]
Our proof of the upper bound on the volume of a periodic orbit in Proposition~\ref{proposition about regulator} applies equally well to a union of periodic orbits of a common discriminant. In other words, the proof of that proposition gives that if $x_1H, \dots, x_nH$ are distinct periodic $H$ orbits with $\disc (x _ 1 H) = \dots = \disc (x _ n H) = D$ then
\begin{equation*}
\sum_ i \vol (x _ i H) \ll D ^ c
\end{equation*}
hence the total volume of all periodic orbits of discriminant $\leq D$ is $\ll D ^{c+1}$.

By the lower bound in Proposition~\ref{proposition about regulator}, the volume of a periodic $H$-orbit is $\gg 1$,\footnote{Of course, this can also be seen directly and elementarily.} and we conclude that $N (D) \ll D ^{c+1}$.
\end{proof}

\begin{proof}[Second proof of upper bound]
For any periodic orbit $\Gamma gH$ with discriminant $\leq D$ and
associated torus $\T$, the point $\iota(\mathfrak{t}) \in V$ defined
by \eqref{iotadef} satisfies $\iota(\mathfrak{t}) = \mathrm{Ad}(g)
\iota(\mathfrak{h})$. By Theorem~\ref{theorem of Tomanov and Weiss},
$g$ may be taken to belong to a fixed compact subset of $G$;
therefore $\iota(\mathfrak{t})$ belongs to a fixed compact subset of
$V$ and has (by definition of discriminant) denominator $\leq D$
with respect to the lattice $V_{\Z}$. It follows the number of
possibilities for $\iota(\mathfrak{t})$ is $\ll D ^{\dim V+1}$.
%
\end{proof}

The lower bound in Proposition \ref{proposition about number of orbits} is closely tied
to the notion of {\em packets}, as well as the implicit actions of adelic groups, which will be discussed further in \cite{ELMV3}.

\begin{proof}[Proof for the lower bound.]
Let $\Gamma g H$ be a periodic orbit, so that $g H g ^{-1} = \T(\R)$
for an anisotropic torus $\mathbf{T} \subset \G$ as in Proposition
\ref{Basic bijection}. The existence of one such orbit is
established in \cite[Thm.~2.13]{Prasad-Raghunathan} .
 Then for $\delta \in \G(\Q)$, the
orbit $\Gamma\delta g H$ is also periodic; moreover, if $\delta_1,
\delta_2$ define distinct classes in $\Gamma \backslash
\G(\Q)/\T(\Q)$, the orbits $\Gamma \delta_1 g H$ and $\Gamma
\delta_2 g H$ are distinct.

It suffices, then, to show that there exist ``many''
elements in $\Gamma \backslash \G(\Q) / \T(\Q)$.  To do this, we shall make some use of adelic language; this is natural, because the simplest way of even showing that $\G(\Q)$ has ``many'' elements
is the fact that $\G(\Q)$ is a lattice in $\G(\adele)$.

In the following argument, $N_0$ will denote a {\em sufficiently large integer}. If it is taken
sufficiently large (depending only on $\G, \Gamma, \T$) then all the statements
of the form ``for $p \geq N_0$'' will be valid.
Let $\adele_f$ be the ring of finite adeles, i.e. $\adele_f$
consists of the restricted product $(x_p) \in \prod \Q_p$ where $x_p
\in \Z_p$ for almost all $p$. Pick an open compact subgroup $K_f =
\prod_{p} K_p$ of $\G(\adele_f)$ containing $\Gamma$. It is a
theorem \cite[Theorem 5.1]{PR} that the ``class number'' $K_f
\backslash \G(\adele_f)/\G(\Q)$ is finite. Pick, therefore, a finite
set of representatives $\omega_1, \dots, \omega_r$ for  $K_f
\backslash \G(\adele_f) / \G(\Q)$. We may assume that  $\omega_j \in
K_p$ for all $p \geq N_0$.

Let $\mathcal{P}$ be a set containing prime numbers larger than
$N_0$; for each $p \in \mathcal{P}$ suppose we are given $g_p
\in\G(\Q_p)$ which does not belong to $K_p \T(\Q_p)$.  For each $p
\in \mathcal{P}$, there are $\delta_p \in \G(\Q), 1 \leq j(p) \leq
r, k_p \in K_f$ so that
$$g_p =  k_p \omega_{j(p)}\delta_p. $$
We claim the associated cosets $\Gamma \delta_p \T(\Q)$ are
disjoint. For, if we had an equality $\Gamma \delta_p \T(\Q) =
\Gamma \delta_q \T(\Q)$, we must in particular, have an equality
$$
 K_f \omega_{j(p)}^{-1}k_p^{-1}g_p \T(\Q) = K_f
 \omega_{j(q)}^{-1}k_q^{-1}g_q
\T(\Q).
$$
 Looking at the ``$p$-component'' of this shows that $g_p
\in K_p \T(\Q_p)$, a contradiction.

%
%

Now let us produce such a collection $\{ g_p \}$.
Let $K$ be a field over which $\T$ splits. Thus $\G$ also splits over $K$.
Let $\alpha$ be a root of $\T$ and $u_{\alpha}: \mathbb{G}_a \rightarrow \G$
the corresponding root subgroup. This morphism of algebraic groups is defined over $K$.
If $p \geq N_0$ is prime and we pick an embedding $K \hookrightarrow \Q_p$ (if one exists),
we obtain in an obvious way an embedding $u_{\alpha}: \mathbb{G}_a \rightarrow \G$
over $\Q_p$.  We note that such an embedding exists for a positive density of $p$,
by the Chebotarev density theorem.

We claim that -- for any such $p$ that is sufficiently large, say $p
\geq N_0$ -- we must have $g_p := u_{\alpha}(p ^{-1}) \notin K_p
\T(\Q_p)$, and moreover the discriminant of the associated periodic
orbit $\Gamma \delta_p g H$ is $\leq c_3 p ^{c_4}$.  These two
claims together complete the proof of the lower bound.

The first claim may be deduced as follows. Suppose that $g_p \in K_p \T(\Q_p)$.  Then there exists
$t_p \in \T(\Q_p)$ so that $u_{\alpha}(p ^{-1}) . t_p ^{-1} \in K_p$.
The image of $u_{\alpha}$ defines a closed subvariety of the affine algebraic variety $\G/\T$. This means
that there exists a regular, right $\T$-invariant function $f_{\alpha}$ on the algebraic variety $\G$,
with the property that it extends the function $u_{\alpha}(x) \mapsto x$.
For sufficiently large $p$ we would necessarily have $f_{\alpha}(K_p) \subset \Z_p$,
contradicting the fact that $u_{\alpha}(p ^{-1}). t_p ^{-1} \in K_p$.

The second claim follows from the following fact. There is a
constant $C > 0$ so that, for all $p \geq N_0$ and for all $K
\hookrightarrow \Q_p$, $\mathrm{Ad}(g_p) . (V_{\Z} \otimes \Z_p)
\subset p ^{-C} (V_{\Z} \otimes \Z_p)$. Here $V_{\Z}$ is as in the
definition of discriminant in Section \ref{discriminantsection}.
Indeed, $$\mathrm{Ad}(g_p) (V_{\Z} \otimes \Z_p) = \mathrm{Ad}\circ
u_{\alpha} (p ^{-1}) (V_{\Z} \otimes \Z_p).$$ Now $\mathrm{Ad} \circ
u_{\alpha}$ defines a morphism of algebraic groups from
$\mathbb{G}_a$ to $\mathrm{GL}(V)$ defined over $K$; choosing a
$K$-basis for $V$, the different matrix entries of $\mathrm{Ad}
\circ u_{\alpha}(x)$ define a collection of one-variable polynomials
in $K[x]$. It would suffice to take for $C$ the largest degree of
any one of these polynomials.
%
\end{proof}




Proposition~\ref{proposition about regulator} and Proposition~\ref{proposition about number of orbits} together imply that the number of periodic $H$-orbits of volume $\leq R$ is at most $\exp (c R)$ for some $c$ (cf. \cite{Oh-finiteness} where the finiteness of the number of such orbits is established).

This bound can be improved in the specific examples; see, for
example, Corollary~\ref{Cor:Explication}, where it is shown that for
$X = \PGL (n, \Z) \backslash \PGL (n, \R)$ the number of periodic
$H$-orbits of volume $\leq R$ is bounded by a bound of the form
$\exp (c R ^{1/2})$ which up to the precise value of $c$ is sharp.

\section {Periodic $H$-orbits and positive entropy}\label{entropy section}

Our main aim in this section is to prove positive entropy for limit
measures arising from periodic $H$-orbits, i.e.\ Linnik's Principle
in Theorem \ref{thm:linnikprinciple}. Throughout this section, let
$\G$ be an $\R$-split algebraic group defined over $\Q$, let $
\Gamma$ be an arithmetic lattice in $G = \G (\R)$ (as always
commensurable with $\G (\Z)$), $H$ be an $\R$-split Cartan subgroup
of $G$, and $\mathfrak h = \Lie H$.

For any $x \in \mathfrak h$ there is attached a one parameter subgroup $a (t) = \exp (t x)$ of $H$, and for Haar measure $\haarG$ on $\Gamma \backslash G$ one can easily show that the metric (or Kolmogorov-Sinai) entropy of the flow is given by
\begin{equation}\label{entropy formula for Haar}
h _ {\haarG} (a ({\cdot})) =  h _ {\haarG} (a (1)) = \sum_ {\alpha \in \Phi} \max (0,\alpha (x))
\end{equation}
where $\Phi$ denotes the set of roots of $G$.

\begin{Theorem}  \label{Entropy for algebraic group}
Let $x \in \Lie H$ and $a (t) = \exp (t x)$ be a
one-parameter subgroup of $H$ as above, and $\rho$ an arbitrary positive real number.
Let $ Y _ i = \left\{ y_{i, 1} H, \dots, y_{i, n _ i} H \right\}$ be a collection of periodic $H$-orbits in $X = \Gamma \backslash G$ and $\Delta _ i \to \infty$ satisfying
\begin{enumerate}
\item \label{discriminant bound condition}
the discriminants of all $y_{i, j} H$ for $j=1,\ldots,n_i$ are at
most $\Delta _ i$
\item \label{volume bound condition}
the total volume of all the orbits in $Y _ i$
is bigger than $\Delta _ i ^ {\rho}$.
\setcounter{savedenumi}{\value{enumi}}
\end{enumerate}
Let $\mu _ i$ be the sum of the volume measures on the periodic $H$-orbits $y_{i,j}H$ in $Y _ i$ ($1 \leq j \leq n _ i$),
 divided by the total volume of these orbits (so that $\mu_i$ is a probability measure).  Suppose that
\begin{enumerate}
\setcounter{enumi}{\value{savedenumi}}
\item \label{escape of mass condition}
$\mu _ i \to \mu$ as $i \to \infty$ in the weak$^*$ topology
for some \emph{probability} measure $\mu$.
\end{enumerate}
Then
\begin{equation} \label{basic entropy bound}
h _ \mu (a ({\cdot})) \geq \rho \min_ {\alpha \in \Phi} \absolute
{\alpha (x)} /2 .\end{equation} If, instead of (\ref{discriminant
bound condition}) above, one assumes
\begin{enumerate}
\item[($\ref{discriminant bound condition}'$)] the discriminant of all $y_{i, j}
H$ for $j=1,\ldots,n_i$ is \emph{equal} to $\Delta _ i$
\end{enumerate}
then \equ{basic entropy bound} can be improved to
\begin{equation*}
h _ \mu (a ({\cdot})) \geq \rho \min_ {\alpha \in \Phi} \absolute {\alpha (x)} \tag{\ref{basic entropy bound}'}
.\end{equation*}
\end{Theorem}

For concreteness, we consider explicitly the case of $G \cong \SL (n, \R)$, $H$ the group of the diagonal matrices, $x \in \mathfrak h$ the diagonal matrix with entries $(n-1)/2, (n-3)/2, \dots,- (n-1) / 2$, and $a (t)$ the corresponding one parameter subgroup of $H$. In this case the roots are $ \alpha (y) = y_i -  y_j$ where $y \in \mathfrak h$ a diagonal matrix with entries $y_1, \dots, y_n$, and hence by \equ{entropy formula for Haar} gives
$h _ {\haarG} (a ({\cdot})) = \binom {n + 1 }{ 3}$. It follows that in the notations of Theorem~\ref{Entropy for algebraic group} we have
\begin{equation*}
h _ \mu (a ({\cdot})) \geq
\begin{cases}
\frac {\rho h_{\haarG} (a ({\cdot}))
}{ 2 \binom {n + 1 }{ 3}} & \text {under assumption (\ref{discriminant bound condition})}\\
\frac {\rho h_{\haarG} (a ({\cdot})) }{ \binom {n + 1 }{ 3}} & \text
{under assumption ($\ref{discriminant bound condition}'$).}
\end{cases}
\end{equation*}
Neither of these bounds seem tight, and indeed in some cases a
better bounds can be obtained. In particular, for $\mathbf G =
\PGL_2$ one can prove under assumption ($1 '$) a sharp estimate $h _
\mu \geq 2 \rho$. Such a bound (which is far from being trivial,
even in this very simple context), or more precisely a $p$-adic
analog of such a bound, can be deduced from Linnik's ``basic lemma''
in \cite{Linnik-book}, and a simplified and explicit derivation of
this bound will be given in \cite{ELMV2}.

We wish to draw attention to assumption (\ref{escape of mass condition}) in Theorem~\ref{Entropy for algebraic group}. If $X = \Gamma \backslash G$ is compact then by passing to a subsequence if necessary this assumption is automatically satisfied. On the other hand, if $\Gamma$ is not a uniform lattice in $G$, whether (\ref{escape of mass condition}) is satisfied or not is a rather interesting issue;
e.g. for $\G = \PGL_n$ it is closely related to analytic properties of Dedekind $\zeta$-functions.

If one is willing to compromise on the quality of the entropy bound
\equ{basic entropy bound} or ($\ref{basic entropy bound} '$), it is
possible to relax slightly the assumption that $\mu (X) = 1$,
assuming only that $\mu (X) > c$ for some explicit $c$ (which for
large $n$ will be extremely close to 1). We do not give the details
as this situation does not seem likely to arise in any natural
context.

\subsection{Entropy and bounds on measures of tubes}

In this and the next subsection we review some well-known facts
about entropy in a form which will be convenient for our purposes
and deduce Theorem \ref{Entropy for algebraic group} from the
following proposition. Note that the fact that we are dealing with
spaces which are not compact causes minor complication in an
otherwise straightforward argument. Moreover, because the space
might not be compact it is necessary to assume that the weak$^*$
limit is a probability measure. 

Let $a(t)$ denote a $\R$-diagonalizable one parameter subgroup of a semisimple group $G= \G (\R)$ as above. Fix some open neighborhood of the identity $B \subset G$. For any $s < t \in \R ^ +$, denote
\begin{equation*}
B ^ {( s,t )} = a (-s) B  a (s) \cap a (- t) B a (t)
.\end{equation*}

\begin{Proposition} \label{entropy proposition}
Suppose $\mu _ i$ is a sequence of $a (t)$-invariant probability measures on $\Gamma \backslash G$ converging in the weak$^ {*}$ topology to a probability measure $\mu$. Suppose further that
there is a sequence of positive real numbers $t _ i \to \infty$ so that, for every compact $\Omega \subset \Gamma \backslash G$, there exists an open neighborhood of the identity $B \subset G$ with:
\begin{equation} \label{equation for which elon proposed a very long name which is not so very long if dictated by voice}
\mu _ i \times \mu _ i \left ( \left\{ (x, y) \in \Omega ^ 2: y \in x B ^ {( -t _ i, t_i )} \right\} \right) < C _ \Omega e ^ {- 2 \eta t _ i}
.\end{equation}
Then the metric entropy of $\mu$ with respect to the flow $a(t)$ satisfies $h _ \mu (a(t)) \geq \eta$.
\end{Proposition}

In particular, we have the following:
\begin{Corollary} \label{entropy corollary}
Let $\mu _ i, \mu $ be $a(t)$-invariant probability measures on
$\Gamma \backslash G$ with $\mu _ i \weaklygoesto \mu$. Suppose that
there is a sequence of positive real numbers $t _ i \to \infty$ so
that, for every compact $\Omega \subset \Gamma \backslash G$, there
exists an open neighborhood of the identity $B \subset G$ with
\begin{equation} \label{pointwise volume bound}
\mu _ i  \left ( x B ^ {( -t _ i, t_i )} \right) < C' _ \Omega e ^ {- 2 \eta t _ i} \qquad \text{for every $x \in \Omega$}
.\end{equation}
Then the metric entropy of the flow $a(t)$ satisfies $h _ \mu (a(t)) \geq \eta$.
\end{Corollary}

\begin{proof} [Proof of Corollary~\ref{entropy corollary} assuming Proposition~\ref{entropy proposition}]
By Fubini,
\begin{equation} \label{reduction step one}
\mu _ i \times \mu _ i \left ( \left\{ (x, y) \in \Omega ^ 2: y \in x B ^ {( -t _ i, t_i )} \right\} \right) = \int_ \Omega \mu _ i \left ( x B ^ {( -t _ i, t_i )}  \cap \Omega\right) d \mu _ i (x)
\end{equation}
applying \equ{pointwise volume bound} we get
\begin{equation*}
\equ{reduction step one} \leq C ' _ \Omega \int_ \Omega e ^ {- 2 \eta t _ i} d \mu _ i (x)\leq C' _ \Omega e ^ {- 2 \eta t _ i}
\end{equation*}
and hence the condition on $\mu _ i$ in Corollary~\ref{entropy corollary} implies that in Proposition~\ref{entropy proposition}.
\end{proof}

As we will see below the assumption in Theorem \ref{Entropy for algebraic group}
regarding the total volume and the information from Proposition \ref{proposition discrete}
regarding the discreteness of periodic orbits are enough to deduce the assumptions to Corollary \ref{entropy corollary}.

\begin{proof} [Proof of Theorem \ref{Entropy for algebraic group} assuming Proposition \ref{entropy proposition}]
The idea in words is that the assumption \eqref{pointwise volume
bound} to Corollary \ref{entropy corollary} expresses an upper bound
on the volume of the intersection of collections of periodic orbits
with small tubes. However, by Proposition~\ref{baby}, pieces of
closed orbits can never be too close to each other, so,
 if the tubes are small enough and the total volume is big enough,
 a ``trivial''  bound on the volume of pieces of periodic orbits is sufficient
 to obtain \eqref{pointwise volume bound} and so positive entropy.

Let $B _ H \subset H$ be a compact neighborhood of the identity and
let $B _ R = \exp V$ where $V$ is a compact neighborhood of zero in
the linear hull of the nonzero root spaces for $\mathfrak{h}$ acting
on $\mathfrak{g}$. Then $B =B _ R B _ H$ is a neighborhood of the
identity $e \in G$. We are going to prove the assumption to
Corollary \ref{entropy corollary} for an arbitrary compact set
$\Gamma \backslash \Gamma \Omega$ with compact $\Omega \subset G$.

For this, note first that
for $t \geq 0$
\begin{equation*}
B ^{(- t, t)} = \bigl(a (t) B _ R a (- t) \cap a (- t) B_R a (t)\bigr) B _ H.
\end{equation*}
Define $\kappa=\min_{\alpha\in\Phi}|\alpha(x)|$. We may assume $a$
is regular and so $\kappa>0$. There exists some constant $c> 0$ such
that
\begin{equation*}
B ^ {(- t, t)} \subset B_{ce ^{-\kappa t}} ^ G (e) B _ H,
\end{equation*}
for all $t > 0$, where $B ^ G_{\delta}(e)$ denotes a $\delta$-neighbourhood of the identity in $G$.

Recall that if $g_1,g_2 \in \Omega B$ give rise to periodic
$H$-orbits of discriminant less than $\Delta_i$ then $d(g_1,g_2) \gg
_ {\Omega}\Delta_i ^{-1}$ unless $g_1 N_G(H)=g_2 N_G(H)$ by
Proposition~\ref{proposition discrete}. Now set $t _ i= \frac { \log
\Delta _ i} {\kappa} +A$ for some constant $A$ that depends on
$B,\Omega$. By choosing $A$ sufficiently large, we can achieve that
$g_1,g_2 \in g B ^{(-t_i,t_i)}$ with
 $\Gamma g_1,\Gamma g_2 \in Y_i$ and $g \in \Omega$ implies $g_1 N_G(H)=g_2 N_G(H)$.

Therefore, $Y_{i}\cap\Gamma g B ^{(-t_i,t_i)}$ is contained in {\em at most one} periodic $N_{G}(H)$-orbit
 $\Gamma g_{i}N_{G}(H)$ and hence in at most $|N_{G}(H)/H|$ periodic $H$-orbits. More precisely,
 there exists $\delta=\delta(\Omega,B)>0$ and $x_{1},\dots,x_{J}\in \Gamma g_{i}N_{G}(H)$ with $J=|N_{G}(H)/H|$ such that
 $$Y_{i}\cap\Gamma g B ^{(-t_i,t_i)}\subset \bigcup_{j=1}^{J}x_{j}B^H_{\delta}(e).$$

Since the total volume of $Y_i$ with respect to a fixed Haar measure of $H$ is assumed to be $\Delta_i ^ \rho$ we see that for the normalized measure
$$
\mu_i\left(\Gamma g B ^{(-t_i,t_i)}\right) \leq \Delta_i
^{-\rho}\sum_{j}\mu_{H,\mathrm{Haar}}(x_{j}B^H_{\delta}(e))\ll
\Delta_i ^{-\rho} \ll e ^{- {\rho \kappa}t_i}.
$$
 This shows \eqref{pointwise volume bound} with $\eta=\frac{\rho\kappa}{2}$ and the theorem in the case
 of \eqref{discriminant bound condition}. For (\ref{discriminant bound
 condition}') the proof uses the first case of Proposition \ref{proposition
 discrete} instead of the second.
\end{proof}

\subsection{Proof of Proposition \ref{entropy proposition}}
For any set $\Omega \subset \Gamma \backslash G$, we let
\begin{equation*}
\Omega _ {M (\epsilon)} = \left\{ x \in \Gamma \backslash G: \sup_ {T \in \Z ^ +} \frac {1 }{ 2T} \sum_ {n = -T+1} ^ {T} 1 _ {\Omega ^ \complement} (x a(-n))  < \epsilon \right\}
\end{equation*}
with $\Omega ^ \complement$ denoting the complement of $\Omega$ in $\Gamma \backslash G$.
Note that for $\epsilon < 1/2$ the above definition particularly implies that $\Omega _ {M (\epsilon )} \subset \Omega$. From the maximal inequality, we know that for any $a (t)$-invariant probability measure $\mu$
\begin{equation*}
\mu (\Omega _ {M (\epsilon)} ^ \complement) \leq 10 \epsilon ^{-1} \mu (\Omega ^ \complement)
.\end{equation*}

For any partition $\mathcal{P}$ of $\Gamma \backslash G$ we let
\begin{equation*}
\mathcal{P} ^ {( s, t )} = \bigvee _ {s \leq n \leq t} \mathcal{P} a (n)
.\end{equation*}
All our partitions will be implicitly assumed to be finite. We will use $[x] _ \mathcal{P}$ to denote the unique element of $\mathcal{P}$ containing the point $x$. For any finite partition $\mathcal{P}$,
we let $H _ \mu (\mathcal{P}) $ denotes its entropy, i.e.
\begin{equation} \label{entropydef}
H _ \mu (\mathcal{P}) = - \sum_ {P \in \mathcal{P}} \mu (P) \log \mu (P)
.\end{equation}

We shall say that the partition $\mathcal{P}$ is $\mu$-regular if for every $P \in \mathcal{P}$ the boundary of $P$ has $\mu$ measure zero.

We now fix a partition $\mathcal{P}$ and a compact set $\Omega \subset \Gamma \backslash G$ which will be used for the remainder of this section. First, we let
$\Omega$ be a compact set so that the $\mu$ measure of the interior of $\Omega$ is $> 1 - \epsilon ^ 2 / 200 $. Then since $\mu _ i \weaklygoesto \mu$, for all sufficiently large $i$
\begin{equation}\label{equation regarding measure of Greek Omega}
\mu _ i (\Omega) > 1 - \epsilon ^ 2 / 100\quad \text{and} \quad
\mu _ i (\Omega _ {M(\epsilon )}) > 1 - \epsilon / 10
.\end{equation}

Let $B _ 1$ be a relatively compact symmetric open neighborhood of the identity in $G$ so that:
\begin{enumerate}
\item $B _ 1 ^ 2 \subset B$ with $B$ a neighborhood of the identity
satisfying our assumption \equ{equation for which elon proposed a
very long name which is not so very long if dictated by voice}
\item for any $ x \in \Omega$,
\begin{equation*}
x. B _ 1 \cap x a (- 1) B _ 1 a (1) = x B _ 1 ^ {( 0, 1 )}
.\end{equation*}
\end{enumerate}
We now take $\mathcal{P}$ to be any (finite) $\mu$-regular partition so that for every $x \in \Omega$ the element $[x]_\mathcal{P} \subset x B_1 $ (note that the existence of such partitions is immediate).

\begin{Lemma} \label{lemma about choice of B}
Without loss of generality we can choose $B _ 1 $ so that for some $C$ and every $n \in \Z ^ +$
\begin{align}
\label{first equation about B}
B_1 ^ {(0, n)} & \subset F_n  \cdot B_1 ^ {( 0, n + 1 )} \qquad \text{for some $ F_n \subset G$ with $\absolute {F_n} \leq C$}\\
\label{second equation about B} B_1 ^ {(-n, 0)} & \subset F_n \cdot
B_1 ^ {( - n - 1, 0 )} \qquad \text{for some $ F_n \subset G$ with
$\absolute {F_n} \leq C$} .\end{align}
\end{Lemma}

\begin{proof}
We will only prove \equ{first equation about B}; the proof of
\equ{second equation about B} is similar. Let $\mathfrak g$ be the
Lie algebra of $G$, and let $\mathfrak g = \bigoplus _{\lambda}
\mathfrak g _ \lambda$ be the decomposition of $\mathfrak g$ into
eigenspaces for $\operatorname{Ad} (a(t))$. The Riemannian metric
defined earlier on $G$ gives us a Euclidean norm on $\mathfrak g$.
So, for any $\delta$ let $B _ \lambda (\delta) = \left\{ \mathbf g
\in \mathfrak g _ \lambda: \norm {\mathbf g} < \delta \right\}$. We
take $B _ 1 = \exp (\sum_{\lambda} B _ \lambda (\delta))$ where
$\delta$ is chosen to be sufficiently small so that $B _ 1 ^ 2
\subset B$ and so that the map $\exp$ is a diffeomorphism from
$\sum_\lambda B_{\lambda}(\delta)$ onto its image. Clearly,
\begin{gather*}
a(t) B _ 1 a(-t) = \exp\Bigl (\sum _ {\lambda} B _ {\lambda} (\exp (t \lambda) \delta) \Bigr) \\
B _ 1 ^ {( s, t )} = \exp \Bigl (\sum _ {\lambda \geq 0} B _ \lambda
(\exp (s \lambda) \delta) + \sum _ {\lambda < 0} B _ \lambda (\exp
(t \lambda) \delta) \Bigr) .\end{gather*} Assuming $\delta$ is
sufficiently small, for every $\mathbf g \in \sum_\lambda B _
\lambda (\delta)$ and $n > 0$ we have:
\begin{align*}
\exp (\mathbf g) & \exp \Bigl (\sum _ {\lambda \geq 0} B _ \lambda (\delta/2) \times \sum _ {\lambda < 0} B _ \lambda (\exp ((n+1) \lambda) \delta/2) \Bigr) \\
& \subset
\exp \Bigl (\mathbf g + \sum _ {\lambda \geq 0} B _ \lambda (\delta) \times \sum _ {\lambda < 0} B _ \lambda (\exp ((n+1) \lambda) \delta) \Bigr);
\end{align*}
This follows easily, e.g. the maps $(X,Y) \mapsto \exp(X+Y)$ and $(X,Y) \mapsto \exp(X) \exp(Y)$
from $\g \times \g$ to $G$ have the same derivative at $0$.

Since clearly $\sum _ {\lambda \geq 0} B _ \lambda ( \delta) \times \sum _ {\lambda < 0} B _ \lambda (\exp (n \lambda) \delta)$ can be covered by a fixed finite number, say $C$, of translates of
\begin{equation*}
\sum _ {\lambda \geq 0} B _ \lambda (\delta/2) \times \sum _ {\lambda < 0} B _ \lambda (\exp ((n+1) \lambda) \delta/2),
\end{equation*}
with $C$ independent of $n$, this lemma follows.
\end{proof}

Using Lemma~\ref{lemma about choice of B}, one easily proves the following using induction:
for integers $k \leq 0 \leq m$ and $x \in \Omega$, if
\begin{equation*}
n = \absolute {\left\{ k \leq i \leq m: x a (- i) \in \Omega ^ \complement \right\}}
\end{equation*}
then $[x] _ {\mathcal{P} {( k, m)}} \subset F_x \cdot B _ 1 ^ {(k,m
)}$ for some $F_x \subset \Gamma \backslash G$ satisfying $\absolute
{F} \leq C ^ n$.
In particular, we have the following:

\begin{Corollary} \label{corollary controlling atoms}
For any
$x \in \Omega _ {M (\epsilon)}$ and $t \in \R ^ +$ there is a finite subset
$F _ {x, t} \subset \Gamma \backslash G $
with $\log \absolute {F _ {x, t}} \ll \epsilon t$ so that
\begin{equation*}
[x] _ {\mathcal{P} ^ {( - t, t )}} \subset F_{x, t} \cdot B_1 ^ {( -t, t)}
.\end{equation*}
\end{Corollary}

\begin{Lemma} \label{main entropy lemma}
For every sufficiently large $i$
\begin{equation*}
H _ {\mu _ i} (\mathcal{P} ^ {( - t _ i, t _ i )}) \geq   \mu _ i (
\Omega _ {M (\epsilon )} ) (2 \eta - c \epsilon) t _ i
+O_{\Omega}(1),
\end{equation*}
with $c$ independent of $\Omega, \epsilon, i$.
\end{Lemma}

\begin{proof}
We will first calculate the entropy of $\mathcal{P} ^ {( - t _ i, t _ i )}$ with respect to the measure $\mu_i' = \mu_i |_ {\Omega _ {M (\epsilon )}}$, i.e. the measure
\begin{equation*} \mu_i'  (A) = \frac {1 }{ \mu_i (\Omega _ {M (\epsilon )})
} \mu_i (A \cap \Omega _ {M (\epsilon )}) .\end{equation*} Let
$\left\{ P _ 1, \dots, P _ N \right\} = \mathcal{P} ^ {( - t _ i, t
_ i )}$, and let $p' _ n = \mu _ i ' (P _ n)$. By
Corollary~\ref{corollary controlling atoms}, for every $P _ n$
intersecting $\Omega _ {M (\epsilon )}$ (in particular every $P _ n$
with $p' _ n > 0$), one has that
\begin{equation*}
P _ n = \bigcup_ {m = 1} ^ {M _ n} Q _ {nm}
\end{equation*}
with the $Q_{nm}$ disjoint, each $Q_{nm} \subset x_{nm} B_1 ^ {( - t
_ i, t _ i )}$ for some $x _ {nm} \in \Gamma \backslash G$, and $M _
n < e ^ {c \epsilon t_i}$. Write $q' _ {nm} = \mu _ i ' (Q _ {nm})$,
and let $\mathcal{Q} = \left\{ Q _ {nm} \right\}$. Then since $\sum_
n p'_n=1$, $\sum_ m q'_{nm}=p'_n$, by convexity of $\log$
\begin{equation*}
H _ {\mu_i'} (\mathcal{P} ^ {( - t _ i, t _ i )}) = - \sum_ n p' _ n \log p' _ n \geq - \log (\sum_ i {p' _ n} ^ 2)
\end{equation*}
and finally
\begin{align*}
\sum_ n {p' _ n} ^ 2 & = \sum_ n \left [ \sum_ m q' _ {nm} \right] ^ 2 \leq \sum_ n  \left [ \sum_ m {q' _ {nm}} ^ 2 \right]  \absolute {M _ n}\\
& \leq e ^ {c \epsilon t_i} \sum_ {nm} { q ' _ {nm}} ^ 2
= e ^ {c \epsilon t_i} \int \mu _ i ' ([x] _ \mathcal{Q}) d \mu _ i ' (x)\\
& \leq e ^ {c \epsilon t_i} \int \mu _ i ' (x B ^ {( - t _ i, t _ i )}) d \mu _ i ' (x) \ll_{\Omega} e ^ {c \epsilon t_i - 2 \eta t _ i}
\end{align*}
by Cauchy-Schwarz, the construction of $\mathcal{Q}$, and the
assumption \eqref{pointwise volume bound}. It follows that
\begin{equation*}
H _ {\mu _ i '} (\mathcal{P} ^ {( - t _ i, t _ i )}) \geq (2 \eta - c \epsilon) t _ i + O_{\Omega}(1)
.\end{equation*}
Let $ p _ n = \mu (P _ n)$.
Then
\begin{align*}
H _ {\mu _ i}   (\mathcal{P} ^ {( - t _ i, t _ i )} ) &= - \sum_ n p _ n \log p _ n
\geq -  \sum_ n p _ n \log p' _ n\\
& \geq - \mu _ i (\Omega _ {M (\epsilon )}) \sum_ n p' _ n \log p' _ n\\
& = \mu _ i (\Omega _ {M (\epsilon )}) H _ {\mu' _ i} (\mathcal{P} ^ {( - t _ i, t _ i )} ) \geq \mu _ i (\Omega _ {M (\epsilon )}) (2 \eta - c \epsilon ) t _ i + O_{\Omega}(1)
.\end{align*}
\end{proof}

\begin{proof} [Proof of Proposition~\ref{entropy proposition}]

By Lemma~\ref{main entropy lemma}, we have, for sufficiently large $i$:
\begin{equation*}
H _ { \mu _ i} (\mathcal{P} ^ {( - t _ i, t _ i )}) \geq
\mu _ i (\Omega _ {M (\epsilon )})(2 \eta - c \epsilon) t _ i  + O_{\Omega}(1) \geq (1 - \epsilon) (2 \eta - c \epsilon) t _ i
\end{equation*}
and using subadditivity of $H _ {\mu _ i}$ and invariance of $\mu _ i$ under $a(t)$ it follows that
for any $n$
\begin{equation*}
H _ { \mu _ i} (\mathcal{P} ^ {( - t _ i, t _ i )}) \leq \sum_ {k = - \floorof {t _ i / n}-1} ^ {\floorof {t _ i / n}}
H _ { \mu _ i} (\mathcal{P} ^ {(- kn, -(k+1)n}) = (2 \floorof {t _ i / n} + 2) H _ { \mu _ i}  (\mathcal{P} ^ {(0, n)}).
\end{equation*}
Since $ \mathcal{P}$ is $\mu$ regular, for any fixed $n$ we have
$H _ { \mu _ i}  (\mathcal{P} ^ {(0, n)}) \to H _ {\mu} (\mathcal{P} ^ {(0, n)}) $.
Thus\begin{align*}
H _ { \mu }  (\mathcal{P} ^ {(0, n)}) & \geq
\limsup_ {i \to \infty} \frac {H _ { \mu _ i} (\mathcal{P} ^ {( - t _ i, t _ i )}) }{ 2 \floorof {t _ i / n} + 2} \\
& \geq (1 - \epsilon) \limsup_ {i \to \infty} \frac {2 (\eta - c
\epsilon) t _ i }{ 2 \floorof {t _ i / n} + 2} = (1 - \epsilon) (
\eta - c \epsilon) n .\end{align*} Dividing both sides by $n$ and
taking the limit as first $ n \to \infty$ and then $\epsilon \to 0$
gives the proposition.
\end{proof}

\section{The action of $H$ on quotients of central simple algebras}\label{dynamics section}

This section discusses further the case when $\G$ arises from the
multiplicative group of an $\R$-split central algebra (e.g., $\G =
\PGL_n$!)

In Section \ref{Section algebras}, we explicate the correspondence of Section \ref{sec:gn},
and discuss (for the only time in this paper) the notion of {\em packet}. This notion
will be further developed in \cite{ELMV3}.

In Section \ref{Subsec:Isolation}, we present and refine some existing theorems about the dynamics of the Cartan action in this setting.
In the setting of $\G = \PGL_3$ this isolation theorem was proved (in different language)
by Cassels and Swinnerton-Dyer.


\subsection{Central simple algebras and $\PGL_n$}  \label{Section algebras}
\subsubsection{Central simple algebras: notation} \label{CSANot}
Let $ D_{\Q}$ be a central simple algebra over $\Q$ of rank $n$,
i.e. $ D_{\Q} \otimes_{\Q} \overline{\Q}$ is isomorphic to
$M_n(\overline{\Q})$ as an algebra. We assume that $D_{\Q}$ is {\em
split over $\R$}, that is to say, $D_{\Q} \otimes_{\Q} \R$ is
isomorphic to $M_n(\R)$.

Let $\G $ be the algebraic group ``$PD_{\Q}^{\times}$'' associated
to the projective group of units in $D_{\Q}$. The Lie algebra $\g$
of $\G$ is identified with the quotient of $D_{\Q}$ by its center,
the bracket operation being $[x,y] := xy - yx$.

Let $\order_D$ be an order inside $ D_{\Q}$, i.e. a $\Z$-module of rank $n ^ 2$ containing $1$ and closed under multiplication. Choose a basis $\{ x_1, x_2, \dots, x_{n ^ 2}\}$ for $\order_D$; the adjoint
action of $D_{\Q}$ w.r.t. this basis yields an embedding $\G \rightarrow \GL_{n ^ 2}$.
The intersection $\G(\Q) \cap \GL_{n ^ 2}(\Z)$ is commensurable with the image of $\order_D ^{\times}$
in $\G(\Q)$; in particular,the image $\Gamma$ of $\order_D ^{\times}$
inside $G = \G(\R)$ is an arithmetic lattice.
For $\g_{\Z}$ we take the image of $\order_D$ in $\mathfrak{g}$,
and finally for $B$ the form
$$B: (X,Y) \rightarrow  \tr(1) \tr(X Y) - \tr(X) \tr(Y),$$ where $\tr(X)$ is the trace of ``left multiplication by $X$ on $ D_{\Q}$''. This defines a bilinear form on $D_{\Q}$ which descends to $\mathfrak{g}$.
With these choices it is easy to check that our earlier assumptions hold true.

Furthermore, we fix a maximal commutative totally real semisimple subalgebra $E_H$ of $D_{\Q} \otimes \R$;
then $E_H$ is isomorphic to $\R ^{n}$.
The centralizer of $E_H$ defines a Cartan subgroup $H \subset G$.

\subsubsection{Parameterization of closed orbits; lower and upper bounds for volumes.}
We shall explicate the parameterization of closed orbits in this particular case. Proofs will be given in
\S \ref{CSAProofs}.

Let $\order$ be a totally real order, i.e.\ a ring which is a finite
free $\Z$-module and so that $\order \otimes \Q$ is a totally real
number field. Let $\sigma_1, \dots, \sigma_n: \order \otimes \Q
\rightarrow \R$ be the different real embeddings.  By means of
$\theta := (\sigma_1, \dots, \sigma_n)$ we may regard $\order$ as a
lattice in $\R ^ n$; by means of
$$\eta: x \mapsto (\log |\sigma_1|, \dots, \log |\sigma_n|)$$ we may map $\order ^{\times}$ (the group of invertible elements in $\order$)
onto a lattice in $\R ^{n-1} = \{(x_1, \dots, x_n) \in \R ^ n: \sum_i x_i \}$.  The latter result is Dirichlet's unit theorem.  We define the discriminant and regulator:
\begin{equation} \label{discregdef} \mathrm{discriminant} \, \order :=\vol(\R ^{n}/\theta(\order))^ 2,   \ \ \mathrm{regulator} \, \order := \vol(\R ^{n-1}/\eta(\order ^{\times}))\end{equation}

\begin{Proposition} [Parameterization of periodic orbits for central simple algebras.] \label{Compact orbits}
In the case of the algebraic group $\G=P D_{\Q}^{\times}
$ defined by a central simple algebra $D_{\Q}$, we have bijections between:
\begin{enumerate}
\item[(1')]  periodic $H$-orbits on $\Gamma \backslash G$.
\item[(2')] $\Gamma$-orbits on pairs $(E, \varphi)$
where $E$ is a subfield of $ D_{\Q}$ of degree $n$, and $\varphi: E
\otimes \R \rightarrow E_H$ is an algebra isomorphism. Here $
\gamma$ maps $(E, \varphi)$ to $(\gamma E \gamma ^ {- 1}, \varphi
\circ \Ad_D(\gamma ^{-1}))$. Here $\Ad_D(\gamma^{-1})$ denotes the
conjugation of elements of $D_{\Q}$ by $\gamma^{-1}$.
\end{enumerate}
The bijection associates to the pair $(E, \varphi)$ the periodic
orbit $\Gamma g H$, where $g \in G$ has the property that the
conjugation $h \rightarrow g h g ^{-1}$ maps $E_H$ to $E$ and
coincides with $\varphi ^{-1}$.

The discriminant of the orbit is a constant multiple of the discriminant of the order $E \cap \order_D$; the volume of the orbit is a constant
multiple of the regulator of $\order_D$.
\end{Proposition}

Using this correspondence, we can obtain versions of some of the
results from Section~\ref{sec:gn} with better exponents.

\begin{Corollary} \label{Cor:Explication} Notations as in the Proposition, we have:  \begin{enumerate}
\item[(i)]
The total volume of all orbits of discriminant $D$ is
$\gg_{\epsilon} D ^{1/2-\varepsilon}$ and the volume of any single
periodic orbit of discriminant $D$ satisfies $V \ll_{\epsilon} D
^{1/2+\varepsilon}$.\footnote{This part of the Corollary will not be proved in this paper, but in \cite{ELMV3}. However, it is not used anywhere else in the present paper.} 
\item[(ii)]
Let $n$ be prime.
The volume $V$ of any periodic orbit of discriminant $D$ satisfies $V \gg (\log D)^{n-1}$.
\end{enumerate}
\end{Corollary}
We expect that the exponents of Corollary~\ref{Cor:Explication}
cannot be improved.

%


To be even more concrete, we specialize further to the case of the
split matrix algebra: take $ D = M_n(\Q)$ and $\order _ D=M_n(\Z)$,
the algebra of $n$ by $n$ matrices. In that case  $\G = \PGL_n, G =
\PGL_n(\R), \Gamma = \PGL_n(\Z)$. Take $E_H$ to be the diagonal
subalgebra, thus $H$ is the diagonal torus.

Let $K$ be a totally real
field of degree $n$.

\begin{Definition} By a {\em lattice} in $K$ we shall simply mean a $\Z$-submodule of $K$ of rank $n$.
Two lattices $L_1, L_2 \subset K$ are called $K$-equivalent (or simply equivalent) if there is $k \in K ^{\times}$ so that $k. L_1 = L_2$.
Attached to any lattice $L$ we have the associated order $\order_L := \{ \lambda \in K: \lambda L \subset L \}$. Then $K$-equivalent lattices have the same associated order.
\end{Definition}


\begin{Corollary} \label{classical bijection}
There is a bijection between:
\begin{enumerate}
\item[(1'')] periodic $H$-orbits on $\PGL_n(\Z)  \backslash \PGL_n(\R)$.
\item[(2'')] Triples $(K, L, \theta)$
of a totally real number field, a $K$-equivalence class of lattices in $K$,
and an algebra
isomorphism $\theta: K \otimes \R \rightarrow \R ^ n$.
\end{enumerate}
Here the triples of (2'') are considered up to isomorphism in the evident sense.

The discriminant of the orbit
is a constant multiple of the discriminant of the order $\order_L$ associated to $L$,
and the volume of the orbit is a constant multiple of the regulator of $\order_L$.
\end{Corollary}
If one identifies $\PGL_n(\Z) \backslash \PGL_n(\R)$
with homothety classes of lattices in $\R ^ n$, the bijection
assigns to $(K,[L], \theta)$ the $H$-orbit of the homothety class
of the lattice $\theta(L)$. Here $L \in [L]$ is arbitrary.

It is often convenient to group the data of [(2'')] by the order $\order_L$. For instance,
in the case when $\order_L$ is a maximal order, the associated equivalence classes $[L]$
are parameterized by the class group of $\order_L$. Thus, approximately speaking,
periodic $H$-orbits on $\PGL_n(\Z) \backslash \PGL_n(\R)$ are parameterized
by an order in a totally real field, together with an ideal class.

\begin{Remark} \label{Packet}
There is a natural weaker equivalence relation on lattices $L \subset K$. Namely,
we say that $L, L'$ are locally equivalent if, for every prime $p$,
there exists $\lambda \in (K \otimes \Q_p)^{\times}$ so that $(L \otimes \Q_p) = \lambda (L' \otimes \Q_p)$.  Via the parameterization above, this groups the periodic orbits
into finite equivalence classes with the following properties:
\begin{enumerate}
\item Periodic orbits in the same equivalence class
have the same volume and same stabilizer in $H$.
\item Each equivalence class has total volume $V$ satisfying
$$D ^{1/2-\epsilon} \ll_{\epsilon} V \ll D ^{1/2+\epsilon}$$
where $D$ is the common discriminant.
\item
The set of periodic orbits in each equivalence class is acted on (simply transitively) by a suitable class group.
\end{enumerate}
This type of grouping can be done (after choosing some auxiliary data) in a much more general setting;
this is the phenomena of {\em packets} that will be discussed further in
\cite{ELMV3}.
\end{Remark}

Again, in this special context we can obtain versions of some
results of Section~\ref{sec:gn} with sharp exponents:
\begin{Corollary} \label{Numberoforbits}
For the $H$-action on $\PGL_n(\Z) \backslash \PGL_n(\R)$:
\begin{enumerate}
\item[(i)]
There is $c > 0$ so that any bounded $H$-orbit
of discriminant $D$ is contained in $\Omega(c D ^{1/2})$; the exponent is sharp.
\item[(ii)]
If $n$ is prime, the number $N(V)$ of periodic $H$-orbits on $\PGL_n(\Z) \backslash \PGL_n(\R)$
of volume $\leq V$ satisfies $$V ^{1/(n-1)} \ll \log N(V) \ll V ^{1/(n-1)}. $$
\end{enumerate}
\end{Corollary}
The behavior of (ii) when $n$ is not prime is heavily influenced by
the existence of fields with intermediate subfields. Since the main
aim of the corollary is simply to contrast with the case of $n=2$,
when one has the asymptotic $\log N(V) \sim V$, we do not attempt to
analyze the general case.


%



\subsubsection{Proofs.} \label{CSAProofs}
Let $\G$ be as in Section \ref{CSANot}.
Then maximal torii in $\G$ are in bijection with degree $n$ subfields of $D_{\Q}$: to each maximal torus $\T \subset \G$ one associates the preimage of its Lie algebra, under the natural map $D_{\Q} \rightarrow \g$.

\begin{proof} (of Prop. \ref{Compact orbits}.)
The first part of the Proposition is a consequence of
Lemma~\ref{Basic bijection}, taking into account the above remark
and the Skolem-Noether theorem \cite[Thm. IX.6.7]{Hungerford} (which
assures that there exists a $g \in G/H$ such that $e \mapsto g ^{-1}
e g$ coincides with $\varphi: E \otimes \R \rightarrow E_H$). For
the assertion concerning discriminant, let $\{ \bar{e}_1, \dots,
\bar{e}_{n-1}\}$ be a $\Z$-basis for the image of $E \cap \order_D$
in the quotient space $E/\Q$; lift them to $e_1, \dots, e_{n-1} \in
E \cap \order_D$. Then $e_0 = 1, e_1, \dots, e_{n-1}$ is a basis for
$E \cap \order_D$. The discriminant $D$ of the periodic orbit
attached to $(E, g H)$ is
$$ \det \{ B(\bar{e}_i, \bar{e}_j) \} = \det
\{ n \tr(e_i e_j) - \tr(e_i) \tr(e_j) \}_{1 \leq i,j \leq n-1}.$$
On the other hand,
$$\disc(E  \cap \order_D) =\det \{\tr(e_i e_j) \}_{0 \leq i,j \leq n-1}  = n  \det \{ \tr(e_i e_j) - \tr(e_i) \tr(e_j)/n \}_{1 \leq i,j \leq n-1}
$$
From this it follows that $D = n ^{n-2} \disc(E \cap \order_D)$.
\end{proof}

\begin{proof} (of Corollary \ref{Cor:Explication}.)
\begin{enumerate}
\item
The proof of (i) is most natural after the introduction of packets. We defer it to a later paper.
It is not necessary for the proof or statement of any other result in this paper.

\item For the proof of (ii), the assertion $V \gg (\log D)^{n-1}$, we prove the corresponding
fact about orders. If $\order$ is a totally real order of discriminant $D$ and regulator $V$ (as in
\eqref{discregdef}), we will show that $V \gg (\log D)^{n-1}$.

Let $\| \cdot \|$ be the sup norm on $\R ^ n$ defined by $\|(x_1,
\dots, x_n)\| = \max_i |x_i|$. We claim that, for any $\lambda \in
\order, \lambda \notin \Z$, we have $\| \theta(\lambda)\| \gg D
^{\frac{1}{n(n-1)}}$. In fact, $\| \theta(\lambda ^ i)\| =  \|
\theta(\lambda)\| ^ i$ and, because $n$ was prime, $1, \lambda,
\dots, \lambda ^{n-1}$ form a $\Z$-basis for the order
$\Z[\lambda]\subseteq\order$. Therefore, the covolume of $\order$ is
at most $C \| \theta(\lambda)\| ^{n(n-1)/2}$, where $C$ depends only
on $n$. Because the covolume of $\order$ is precisely $D ^{1/2}$,
the claim follows.

Therefore, every nonzero vector in $\eta(\order ^{\times})$ has length $\gg_n \log D$.
Now $\eta(\order ^{\times})$ has rank $n-1$, so it must have covolume $\gg_n (\log D)^{n-1}$, as required.
\end{enumerate}
\end{proof}

\begin{proof} (of Corollary \ref{classical bijection}).
It follows from Proposition \ref{Compact orbits}: periodic $H$-orbits
are in bijection with $\Gamma$-orbits on pairs $(K, \theta)$, where $K \subset M_n(\Q)$
is a totally real field. Given such a pair, we associate to it the triple
$(K, [L], \theta)$.  Here $L$ is obtained in the following way: fix a nonzero $a \in \Q ^ n$ and let $\jmath(x)= x a$ for $x \in K$.  Now set $L = \jmath ^{-1} (\Z ^ n)$.

Conversely, given a triple $(K, [L], \theta)$,
choose a basis $\lambda_1, \dots, \lambda_n$ for $L$ to obtain an embedding
$\iota: K \rightarrow M_n(\Q)$ with the property that $$K \cap M_n(\Z) = \order_L =
\{ \lambda \in L: \lambda L \subset L \}.$$

We thereby obtain a pair $(\iota(K), \theta)$ as in the previous section, and moreover
the $\Gamma$-orbit of this pair is independent of the choice of basis.
\end{proof}
\begin{proof} (of Corollary \ref{Numberoforbits})
(i) Let $(K, [L], \theta)$ be data parameterizing a periodic orbit;
and let $L\subset K$ be any lattice in $K$ in the class $[L]$.
%
If $\alpha \in L$, then also $\alpha  \order_L \subset L$; here
$\order_L$ is the order associated to $L$, as in
Corollary~\ref{classical bijection}. Therefore $\theta(\alpha)
.\theta(\order_L) \subset \theta(L)$. Let $D$ be the discriminant of
$\order_L$, and $\nu$ the volume of $\R ^ n/\theta(L)$. Computing
covolumes, this shows for any $(x_1, \dots, x_n) \in \theta(L)$, we
have $\prod_{i} |x_i| D ^{1/2} \geq \nu$.

Writing $\|(x_1, \dots, x_n)\| = \max_i |x_i|$, we can rephrase in the following way: for any $\mathbf{x} = (x_1, \dots, x_n) \in H. \theta(L)$, $ \| \mathbf{x}\| ^ n D ^{1/2} \geq \nu$.
It follows from this and lattice reduction that any lattice $\theta(L') \in H. \theta(L)$ has a reduced basis $v_1, \dots, v_n$
with $$\nu ^{1/n} D ^{-1/(2n)} \ll \| v_1 \| \ll \| v_2 \| \ll \dots \ll \| v_n \|;   \ \ \prod_i \| v_i \| \asymp \nu$$
where the implicit constants in $\ll$ depend only on $n$.  But this means also that
$\| v_n \| \ll \nu ^{1/n} D ^{(n-1)/(2n)}$.

Now let $\Gamma g H$ belong to the associated periodic orbit. Then $g ^{-1} M_n(\Z) g$ consists
of endomorphisms $X \in M_n(\R)$ that preserve $\theta(L')$, for some $L' \in H. L$.
If $X \in M_n(\R)$ preserves $\theta(L')$ and is nonzero, there exists $1 \leq i \leq n$ so that
$v_i X  \neq 0$. But $\| v_i \| \ll \nu ^{1/n} D ^{(n-1)/(2n)}$ and $\| v_i X \| \gg \nu ^{1/n} D ^{-1/(2n)}$.
So the operator norm of $X$ w.r.t. $\| \cdot \|$ is $\gg D ^{-1/2}$. This implies that
$\Gamma g H \in \Omega(c D ^{1/2})$ for suitable $c = c(n) >0$, in the notation of \eqref{omegadef}.

(ii)
Let $n$ be prime; we will show that the number $N(V)$ of periodic orbits
of volume $\leq V$ satisfies $V ^{1/(n-1)} \ll \log N(V) \ll V ^{1/(n-1)}$.

The upper bound follows from the bound $V \gg (\log D)^{n-1}$, established
in Corollary \ref{Cor:Explication}, as well as
Proposition
\ref{proposition about number of orbits}.

The lower bound follows from Corollary \ref{classical bijection}
and the following fact,  also used in Section \ref{examples section}:
There exists $a_n > 0$ and $\gg X ^{a_n}$ totally real fields $K$
with $\disc(K) \leq X$, each has a maximal order $\order_K$ satisfying
$$\mathrm{regulator}(\order_K) \ll (\log \mathrm{discriminant}(\order_K))^{n-1}$$
and the implicit constants depends only on $n$.  See  \cite[\S 3, Prop. 1] {Duke}.
\end{proof}

\subsection{Measure classification and isolation theorems for central simple algebras.} \label{Subsec:Isolation}
Let $G = \mathbf G (\R)$ be a real algebraic group, $\Gamma < G$ a lattice and $H$ a maximal $\R$-split torus. The dynamics of $H$ is drastically different if $\dim H = 1$ or $\dim H > 1$, and this is reflected by the behavior of the periodic $H$-orbits.

We study in detail the case $ \mathbf G = \SL (2)$ in \S\ref{rank one section}. In this case, and in all cases $\dim H = 1$, the behavior of \emph {individual} periodic $H$-orbits in $\Gamma \backslash G$ is quite arbitrary.
When $\dim H \geq 2$ the situation changes dramatically. First consider the case of $\PGL_n(\Z) \backslash \PGL_n(\R)$.
Conjecturally, the following is expected to hold:


\begin{Conjecture} \label{Margulis-conjectures}
Let $H$ be the diagonal split Cartan subgroup in $\PGL_n(\R)$ for $n \geq 3$
and consider its right action on $\PGL_n(\Z) \backslash \PGL_n(\R)$.
\begin{enumerate}
\item Every ergodic $H$-invariant probability measure is algebraic, i.e.~coincides with the ($L$-invariant) volume measure of a closed $L$-orbit for some subgroup $H \leq L \leq \PGL_n( \R)$. \label{measure theoretic part of conjecture}
\item Every bounded $H$-orbit is periodic. \label{topological part of conjecture}
\item For every compact $\Omega \subset \PGL_n(\Z) \backslash \PGL_n(\R)$ there are only finitely many periodic orbits that are contained in $\Omega$. \label{periodic orbits part of conjecture}
\end{enumerate}
\end{Conjecture}

Part (\ref{measure theoretic part of conjecture}) of this conjecture is a special case of more general conjectures by Furstenberg (unpublished), Katok and Spatzier \cite{Katok-Spatzier}, and Margulis \cite{Margulis-conjectures}, and in particular follows from \cite[Conj. 2]{Margulis-conjectures}(properly interpreted); see \cite{Einsiedler-Lindenstrauss-ICM} for more details and more general conjectures).

Part (\ref{topological part of conjecture}) of this conjecture can be traced back to Cassels and Swinnerton-Dyer \cite{Cassels-Swinnerton-Dyer}, but in the form given here is due to Margulis (see e.g. \cite{Margulis-Oppenheim-conjecture}). It would follow from part (\ref{measure theoretic part of conjecture}); this can be shown quite readily using the techniques of \cite{Cassels-Swinnerton-Dyer}. We show below how to deduce part (\ref{topological part of conjecture}) from part (\ref{measure theoretic part of conjecture}) using a more refined and recent result.

The same techniques also shows that part (\ref{periodic orbits part of conjecture}) of the conjecture follows from part (\ref{topological part of conjecture}). This part of the conjecture has also been highlighted by Margulis (see e.g. \cite [Prob. 30]{ Open-problems}).


One of the important results proved in \cite{Cassels-Swinnerton-Dyer} is an isolation theorem: given any periodic $H$ orbit in $X = \PGL_n(\Z) \backslash \PGL_n(\R)$, and any compact subset $\Omega \subset X$, for any $x \in X$ sufficiently close but not on this periodic $H$ orbit, the orbit $xH$ will intersect the complement of $\Omega$. 
This has been strengthened by E. L. and B. Weiss \cite{Lindenstrauss-Barak} for $\PGL_n(\Z)$ and generalized by G. Tomanov to the split central simple algebra case \cite{Tomanov}. 

Here will need the following slight
variant of these results:

\begin{Theorem} \label{Lindenstrauss-Barak}
Let notation be as in Section \ref{CSANot}, so that $\Gamma$ is a lattice in $G = \PGL_n(\R)$
associated to a central simple algebra.
Let $Y \subseteq X := \Gamma \backslash G$ be a closed $H$-invariant set.  Assume that $Y$ contains a closed orbit $yL$ for a reductive group $L \leq G$ such that
\begin{enumerate}
\item $H \leq L$,
\item $yL$ has a finite $L$-invariant volume, and
\item $\overline{ Y \setminus (yL)}$ contains $yL$. \label{non isolation assumption}
\end{enumerate}
Then there exists a strictly bigger reductive $M > L$ such that $y M $ is closed, has finite $M$-invariant volume, and $yM \subseteq Y$. Moreover, if $D_{\Q}= M_n (\Q)$ so that $\G = \PGL_n$, then $Y$ cannot be compact.
\end{Theorem}

The difference between \cite{Lindenstrauss-Barak, Tomanov} and Theorem~\ref{Lindenstrauss-Barak} is that in \cite{Lindenstrauss-Barak, Tomanov} the set $Y$ is assumed to be of the form $Y = \overline {x H}$ for some $x \in \Gamma \backslash G$ (in which case assumption (\ref{non isolation assumption}) is automatic unless $Y=yL$).

We also remark that, if $n$ is prime, the only possible $L$ as in the statement of the Theorem
are $L = H$ and $L=G$.  The preimage of the Lie algebra of $L$ in $D_{\Q}$
is actually a semisimple algebra (this statement may be verified over $\R$; now, note
that all reductive subgroups of $\PGL_n(\R)$ that contain a Cartan subgroup are, in fact, Levi subgroups).  So the preimage of the Lie algebra of $y L y ^{-1}$ in $D_{\Q}$
is a semisimple subalgebra $\mathfrak{l} \subset D_{\Q}$.  In order that $y . L$
have finite invariant measure, the group $L$ should not have any $\Q$-characters; this means
that $\mathfrak{l}$ must be be a {\em simple} $\Q$-algebra; it has dimension
$e ^ 2$ over its center, some number field $E$. Then $\mathfrak{l} \otimes \C$
is a sum of $[E:\Q]$ copies of $M_e(\C)$; every irreducible representation of it
has dimension $e$, so $e|n ^ 2$. Since $n$ is prime, this forces $e=1$ and $\mathfrak{l} = E$,
or $e=n$ and $\mathfrak{l} = D_{\Q}$. In the former case, $L=H$; in the latter case, $L=G$.

\begin{proof}  [Proof of Theorem \ref{Lindenstrauss-Barak}]
We identify $G \cong \PGL_n(\R)$ and take $H$ to be the group of diagonal $n \times n$ matrices (with proportional matrices identified). For $a _ 1, \ldots, a _ n$ we let $\operatorname{diag} (a _ 1, \ldots, a _ n)$ denote the diagonal matrix with entries $a _ 1, \ldots, a _ n$.

The proof of \cite[Thm.~1.1]{Lindenstrauss-Barak} given in \cite [pp. 1490--1492] {Lindenstrauss-Barak} actually proves Theorem~\ref{Lindenstrauss-Barak}, not just for central simple algebras but for any $\Gamma < \PGL_n(\R)$, \emph {provided the following conditions are satisfied:}

\begin{itemize}
\item
The set $Y$ contains a periodic $H$-orbit
$\Gamma g _ 0 H$;
\item the orbit of $\Gamma g _ 0 $ with respect to each of the $\R$-subtori
\begin{equation*}
H_{ij}=\{ \operatorname{diag} (a _ 1, \ldots, a _ n) \in H: a _ i = a _ j \}
\end{equation*}
is not closed.
\end{itemize}

By \cite[Thm.~2.13]{Prasad-Raghunathan} the orbit $yL$ contains a periodic $H$-orbit $ \Gamma g_0 H$.
(In \cite{Prasad-Raghunathan} the group is actually assumed to be semisimple. However,
$L$ is the product of its center and the semisimple commutator and the same holds up to finite index
for its arithmetic lattice. Therefore, the center of $L$ must have compact orbit by the finite volume assumption.)

By Proposition \ref{Compact orbits} and the assumption that $H$ is $\R$-split the periodic orbit $\Gamma g _ 0 H$ corresponds to a totally real number field $K \subset D$ of degree $n$. Now consider two indices $1 \leq i < j \leq n$ and the $ \R$-subtorus $H _{ ij}$ as above.  Then $\Gamma \cap g _ 0 H _{ij} g _ 0 ^ { - 1} = \{ \xi \in K \cap \order _ D ^ \times: \phi(\xi) _ i = \phi(\xi) _ j \} $ where $\phi:K \rightarrow \R ^ n$ is the algebra homomorphism induced by conjugation by $g _ 0$
as in Lemma \ref{Compact orbits}. Since $\phi (\xi) _ k$ for $k = 1,\ldots, n$ are precisely the various Galois embeddings of $\xi$ into $\R$ and at least two of these give the same value, $\xi$ cannot generate $K$ and so for any $k=1,\ldots, n$ there is a second index $\ell$
with $\phi (\xi) _ k = \phi (\xi) _\ell$.  In particular,
$\Gamma \cap g _ 0 H _{ij} g _ 0 ^ { - 1}$ cannot be a lattice in
$g _ 0 H _{ij} g _ 0 ^ { - 1}$ since $n \geq 3$.  This establishes the second
assumption above.
\end{proof}

Using Theorem~\ref{Lindenstrauss-Barak} it is easy to prove that (\ref{measure theoretic part of conjecture}) implies (\ref{topological part of conjecture}) and (\ref{topological part of conjecture}) implies (\ref{periodic orbits part of conjecture}) in Conjecture \ref{Margulis-conjectures}.

Indeed,suppose (\ref{measure theoretic part of conjecture}) in  Conjecture \ref{Margulis-conjectures} holds and let $\PGL_n(\Z) g_0 H$ be a non-periodic orbit with compact closure $Y$. Then the restriction of the $H$-action on $Y$ has an $H$-invariant and ergodic probability measure $\mu$ which has to be algebraic by our assumption, i.e.~$\mu$ is the finite $L$-invariant volume on a closed $L$-orbit for some group $L<G$. Since $H \leq L$ it follows that $L$ is reductive. By Theorem \ref{Lindenstrauss-Barak} it follows that $Y$ cannot be compact --- a contradiction.

Suppose (\ref{topological part of conjecture}) holds and there are contrary to (\ref{periodic orbits part of conjecture}) infinitely many periodic $H$-orbits within a particular compact set $\Omega$. Let $Y$ be the closure of the union of a sequence of different periodic orbits $\Gamma g_i H$ in $\Omega$. Since $\Omega$ is compact there exists a limit point $x$ to a sequence $x_i \in \Gamma g_i H$. Clearly $xH \subset \Omega$. By our assumption $x$ has itself a compact orbit. Therefore, the assumption to Theorem \ref{Lindenstrauss-Barak} is  satisfied with $L=H$ which gives a contradiction to $Y$ being a compact invariant set.

\medskip

The best known results towards Conjecture~\ref{Margulis-conjectures}.(\ref{measure theoretic part of conjecture}) are under an additional assumption, namely  positive entropy. Also note that this conjecture is related to another open question by Furstenberg: What are the probability measures on the circle group $\R/\Z$ invariant under multiplication by $2$ and $3$? The best known result there is Rudolph's theorem \cite{Rudolph-2-and-3} saying that an ergodic such measure with positive entropy must be the Lebesgue measure. Analogues in higher dimensions and homogeneous spaces were first obtained by Katok and Spatzier \cite{Katok-Spatzier, Katok-Spatzier-corrections} under additional assumptions, see also \cite{Kalinin-Spatzier}.
More recently, two methods were developed \cite{Einsiedler-Katok, Einsiedler-Katok-II}, \cite{Lindenstrauss-Quantum} for the homogeneous case and together they were sufficient to show that for the $\R$-split Cartan action on $\SL(n,\Z) \backslash \SL(n,\R)$  positive entropy for some element is sufficient to deduce that the measure is the volume measure of the homogeneous space \cite [Cor.~1.4]{ Einsiedler-Katok-Lindenstrauss}.

For our purposes, we would like to consider more generally any lattice attached to an order in a central simple algebras as above and not just the special case of $\SL (n, \Z)$. This extension (which was suggested to us by Silberman and Tomanov) does not pose additional technical difficulties, and we present it below. Note that examples due to M. Rees \cite{Rees-example} show that the theorem below as stated is false for a general lattice $\Gamma$ even if $G = \PGL_3( \R)$.

\begin{Theorem} \label{measure rigidity for central simple algebra}
Let notation be as in Section \ref{CSANot}, so that $\Gamma$ is a lattice in $G = \PGL_n(\R)$
associated to a central simple algebra. Suppose $\mu$ is an $H$-invariant and ergodic probability measure on $X=\Gamma \backslash G$ such that the metric entropy $h_\mu(a)$ with respect to $\mu$ is positive for some element $a \in H$. Then there exists a reductive $L \leq G$ such that $\mu$ is the $L$-invariant volume form on a single periodic $L$-orbit.
\end{Theorem}

We are following the scheme of proof of \cite [Thm.~1.3]{ Einsiedler-Katok-Lindenstrauss} which uses as the main tool the more general result \cite[Thm.~2.1]{ Einsiedler-Katok-Lindenstrauss} about the structure of conditional measures.  We refer to \cite{Lindenstrauss-Quantum} or \cite [Sec.~2.1]{ Einsiedler-Katok-Lindenstrauss} for the basic theory of conditional measures on foliations.  We only recall that for every root of $\PGL _ n (\R)$, or more concretely every pair of indices $1 \leq i,j \leq n$ with $i \ne j$, there exists a system of conditional measures $\mu _ x ^{ij}$ for almost every $x \in X$ defined on the corresponding unipotent subgroup $U _{ij}$.  Positive entropy means precisely that one such system is nontrivial, i.e.\ for some pair $i,j$ the conditional measures $\mu _ x ^{ij}$ are not just supported on the identity of $U _{ij}$.  If the conditional measures $\mu _ x ^{ij}$ equal the Haar measure of $U _{ij}$ (again for almost every $x$), then $\mu$ is in fact invariant under $U _{ij}$.  Since the dynamics of the unipotent $U _{ij}$ is much better understood, this situation is desirable.  Towards that \cite[Thm.~2.1]{ Einsiedler-Katok-Lindenstrauss} says that for any $i,j$ one of the following three possibilities take place:
\begin{enumerate}
\item[(i)]  The conditional measures $\mu _ {x}^{ij}$ and $\mu_{x}^{ji}$ are trivial
a.e.
\item [(ii)] The conditional measures $\mu _ {x}^{ij}$ and $\mu_{x}^{ji}$ are
Haar a.e., and $\mu$ is invariant under left multiplication with elements of
$L_{ij}=\langle U_{ij},U_{ji}\rangle$.
\item [(iii)] \label{exceptional returns} Let
$A _ {ij}' = \left\{ \operatorname{diag} (s_1,\ldots,s_n): \mathbf s \in (\R ^ \times) ^ n \text{ with }
s_i=s_j \right\}$. Then a.e.\ ergodic component of $\mu$ with respect to $A' _
{ij} $ is supported on a single $C(L_{ij})$-orbit, where $C(L_{ij})=\{ g:
gh=hg$ for all $h \in L_{ij}\}$ is the centralizer of $L_{ij}$.
\end{enumerate}

\begin{Lemma} \label{No exceptional returns}
For $X = \Gamma \backslash \PGL _ n (\R)$
as in Theorem \ref{measure rigidity for central simple algebra} case (iii) of the above
is impossible.
\end{Lemma}

\begin{proof}
In \cite [Thm.~5.1]{Einsiedler-Katok-Lindenstrauss} it has been shown
for $G=\SL _ n(\R)$
that in case (iii) there exists an element $\gamma \in \Gamma$ with the following properties:
\begin{enumerate}
\item diagonalizable over $\R$,
\item $\pm 1$ is not an eigenvalue of $\gamma$, and
\item all eigenvalues of $\gamma$ are simple except precisely one which has multiplicity two.
\end{enumerate}
Note that by rescaling its elements a lattice in $\PGL_n(\R)$ gives rise to a lattice in $\SL_n(\R)$. The resulting quotients are isomorphic unless $n$ is even and $\Gamma$ contains elements of negative determinant in which case we get a double cover. Because of this we can use the above result also for $G=\PGL _ n(\R)$, for this we define the eigenvalue of $\gamma \in \PGL _ n(\R)$ as the eigenvalue of the matrix after normalizing the determinant to be $\pm 1$.

Since $D_{\Q} \otimes \R \cong M_n (\R)$ the eigenvalues of left multiplication by $\gamma \in \order_D ^ \times$ on $D_{\Q}$ are precisely the eigenvalues of $\gamma$ when considered as a matrix in $\GL _ n (\R)$, and the multiplicity on $D $ is precisely $n$ times the multiplicity in the matrix. The characteristic polynomial $p (t)$ of left multiplication by $\gamma$ on $D$ therefore factorizes as $p(t)=(t-\xi) ^{2 n} \prod_ {i=3} ^ n (t - \xi _ i) ^ n$  where $\xi, \xi _ 3,\ldots, \xi _ n$ are the pairwise different eigenvalues of $\gamma$.  Since $\gamma \in \order ^ \times$, the polynomial $p (t)$ has integer coefficients and trailing coefficient $\pm 1$.  If $\xi \not \in \Q $, then a Galois conjugate of $\xi$ must be a root of $p(t)$ of the same multiplicity which is impossible by (3).
Therefore, $\xi \in \Q$ which forces $\xi =\pm 1$ --- a contradiction to (2).
\end{proof}

\begin{proof}  [Proof of Theorem \ref{measure rigidity for central simple algebra}]
As explained before our assumption of positive entropy translates to the nontriviality
of some conditional measures, i.e.\ there are pairs of indices $i,j$ for which (i) above fails.
By the above lemma (iii) can never hold, so $\mu$ is invariant under $L_{ij}$.
Let $L$ be the subgroup generated by all the unipotent subgroups $U_{ij}$
under which $\mu$ is invariant. Clearly $L$ is normalized by $H$, and in fact
we can reorder the indices such that $L$ equals $\prod_{i =1} ^ \ell \SL _{m_i} (\R)$ embedded into $\PGL _ n (\R)$ as block matrices.
We will show below that $\sum_ {i=1} ^ \ell m_i=n$.

By \cite [Theorems (a) and (b)] {Margulis-Tomanov-2}, applied to $\mu$
and the group $HL$, we know that there is some $\tilde L \geq L$ which is normalized by $H$ so that
almost every $L$-ergodic component of $\mu$ is the $\tilde L$-invariant measure on a
closed $\tilde L$ orbit. In particular $\mu$ is $\tilde L$-invariant, which unless $\tilde L \leq HL$
contradicts the definition of $L$ since otherwise $\tilde L$ contains further unipotent
subgroups which preserve $\mu$. Let now $x=\Gamma g$ have a closed $\tilde L$-orbit $x \tilde L$ of finite volume. Then $\Lambda_{\tilde L}=g \tilde Lg ^{-1}\cap \Gamma$ is a lattice in $g \tilde Lg ^{-1}$, and
so the latter is defined over $\Q$. Therefore, the same is true for the
semi-simple $gLg ^{-1}=[g \tilde Lg ^{-1},g \tilde Lg ^{-1}]$, $\Lambda_L=gLg ^{-1}\cap \Gamma$
is a lattice in $gLg ^{-1}$, and $xL$ is closed with finite volume. However,
this implies $\tilde L=L$.

Moreover, by \cite [Theorems (a) and (b)] {Margulis-Tomanov-2} $\mu$
is supported by a single orbit of the normalizer $N_G(L)$.
If in addition  $\sum_ {i=1} ^ \ell m_i=n$, then $N_G(L)=HL$ and $\mu$ must be the unique $HL$-invariant measure on the $HL$-orbit. Since $HL$ is reductive, this would prove the theorem.

So suppose $\sum_ {i=1} ^ \ell m_i<n$. Since $gLg ^{-1}$ is a $\Q$-group, so are the normalizer $N_G(gLg ^{-1})$ and its center $Z=C(N_G(gLg ^{-1}))$. Here $g ^{-1}Zg<H$ consists of all diagonal matrices
for which all entries in a block corresponding to one of the factors of $L$ are equal to all other entries in the same block, and for which all remaining entries
are also equal. Moreover, if $\mathbf{L}', \mathbf{Z}<\D ^ \times$ denotes the algebraic group over $\Q$ corresponding to $gLg ^{-1}$ and $Z$, then the Lie algebra $\mathfrak v$ of $\mathbf{L'}\mathbf{Z}$ is a $\Q$-subspace of $D_{\Q}$. Going back to $\GL_n(\R)$ and the block matrix description of all of the above groups, it is easy to see that $\mathfrak v$
is invariant under left multiplication by $\mathbf{Z} \subset \D ^ \times$.  The determinant of this representation defines a $\Q$-character of $\mathbf{Z}$.
We claim that this determinant is as a character linearly independent from the norm character (defined by the left multiplication on $D_{\Q}$ as a representation). However, this again can be easily seen from the block matrix description: the entries corresponding to the factors of $L$ are used in higher powers than the remaining ones. Therefore, $\mathbf{Z}$ has $\Q$-rank at least two --- modulo $\z ^ \times$ this shows that $Z$ has $\Q$-rank at least one. If $h \in g ^{-1}Z g$ is such that the value of the $\Q$-character is bigger than one, then it follows that $h ^ j.x \rightarrow \infty$ for $j \rightarrow \infty$ for all points in the $N_G(L)$-orbit (where the divergence to infinity is understood within the $N_G(L)$-orbit). This contradicts Poincar\'e recurrence,  and so proves the theorem.
\end{proof}

\section{Density and distribution of periodic orbits in higher rank}\label{applications section}

In this section we employ the tools developed in \S\ref{entropy section} and \S\ref{dynamics section} to prove statements about periodic orbits in $\Gamma \backslash G$ with $\Gamma$ a lattice associated to an order in a degree $n$ central simple division algebra, split over $\R$, in $G \cong \PGL_n(\R)$.

We first use Theorem~\ref{Entropy for algebraic group} to give a general, if somewhat messy, theorem regarding the distribution of periodic $H$-orbits. Note that if $\Gamma \backslash G$ is noncompact, it is also a conditional result.
Later in this section we will deduced from this theorem somewhat cleaner results regarding density properties of these periodic orbits.

\begin{Theorem} \label{comprehensive equidistribution theorem}
Let $\Gamma$ be a lattice in $G \cong \PGL_n(\R)$ obtained from an order in a central simple algebra of degree $n$, split over $\R$. Let  $\rho>0$ be arbitrary.
Let $ Y _ i = \left\{ y_{i 1} H, \dots, y_{i n _ i} H \right\}$ be a collection of periodic $H$-orbits in $X = \Gamma \backslash G$ and $\Delta _ i \to \infty$ satisfying
\begin{enumerate}
\item
the discriminant of all $y_{i j} H$ is
at most $\Delta _ i$
\item
the total volume of all the orbits in $Y _ i$
is bigger than $\Delta _ i ^ {\rho}$.
\setcounter{savedenumi}{\value{enumi}}
\end{enumerate}
Let $\mu _ i$ be the sum of the volume measures on the compact $H$-orbits $y_{ij}H$ in $Y _ i$ ($1 \leq j \leq n _ i$), divided by the total volume of these orbits (so that $\mu_i$ is a probability measure).  Suppose that
\begin{enumerate}
\setcounter{enumi}{\value{savedenumi}}
\item
$\mu _ i \to \mu$ as $i \to \infty$ in the weak$^*$ topology
for some \emph{probability} measure $\mu$.
\end{enumerate}
Then
\begin{equation*}
\mu = \sum_ {i = 0} ^ m a_i \nu_i
\end{equation*}
with $m \in \N \cup \left\{ \infty \right\}$ and
\begin{enumerate}
\item[(a)] each $\nu _ i$ is a $H$-invariant probability measure on $X$,
and $a _ i \geq 0$.
\item[(b)] for every $i \geq 1$ the measure $\nu _ i$ is the $L_i$ invariant probability measure on a single periodic orbit of some closed subgroup $L_i \leq G$ which properly contains $H$.
\item[(c)] $\sum_ {i=1} ^ m a_i \geq c_n \rho$, where $c _ n$ is a constant which depends only on $n$.
\end{enumerate}
Explicitly, we can take\footnote{If one assumes the discriminant of all $y_{i j} H$ is
precisely $\Delta _ i$, this constant can be improved by a factor of two. In either case this bound is far from optimal -- and indeed a better bound can be obtained if one makes a more careful analysis.}
\begin{equation*}
c _ n = \frac {1 }{ 2 \binom {n+1 }{ 3}}
.\end{equation*}
In particular, $\mu$ is not compactly supported.
\end{Theorem}

\begin{proof}
Decompose $\mu$ into its $H$ ergodic components,
\begin{equation*}
\mu =  \int \nu d \tau (\nu)
\end{equation*}
where $\tau$ is a probability measure on the space of (Borel) $H$-invariant and ergodic probability measures on $X$. Then, if $a(t)$ is a one-parameter subgroup of $H$,
\begin{equation*}
h _ \mu (a ({\cdot})) = \int h _ {\nu} (a ({\cdot })) d \tau (\nu)
.\end{equation*} By Theorem~\ref{measure rigidity for central simple
algebra}, each such $\nu$ with $h_{\nu}(a({\cdot}))> 0$ is
algebraic: i.e. is the $L$-invariant probability measure on a single
periodic $L$ orbit, with $L$ a closed group properly containing $H$
(depending on $\nu$), and furthermore each such $\nu$ is not
compactly supported. Since there are only finitely many
possibilities for $L$, and any such $L$ has only countably many
periodic orbits\footnote {This countability issue is mostly
irrelevant for our purposes, but for completeness suppose  $xL$ is
periodic for $x = g \Gamma$. Either by \cite [Thm
1.1]{Ratner-Annals} (where one needs first to show that $x$ is
periodic also under $L'=[L,L]$) or, more directly, by \cite[Prop.
2.1]{ DM} there are only countably many possibilities for $gLg
^{-1}$. But $\left\{ g': g' L g'^{-1} = gLg ^{-1} \right\} = g
N_G(L)$ and $[N_G(L):L]< \infty$. Hence given a conjugate $\tilde L$
of $L$ there are only finitely many orbits $xL$ in $X$ for which $x$
can be written as $\Gamma g L$ with $g L g ^{-1} = \tilde L$.} it
follows that we can write $\mu$ as
\begin{equation*}
\mu = a_0 \nu _ 0 + \sum_ i a_i \nu _ i
\end{equation*}
with $a_i \geq 0$, \ $h _ {\nu _ 0} (a ({\cdot})) = 0$ and each $\nu _ i$ an $L_i$ invariant probability measure on a single periodic $L_i$ orbit, $L_i >H$ a closed subgroup of $G$. This establishes (1) and (2) in Theorem~\ref{comprehensive equidistribution theorem}.

By Theorem~\ref{Entropy for algebraic group}, we know that $h _ \mu (a ({\cdot})) \geq \rho c_n h _ {\haarG} (a ({\cdot}))$ for
\begin{equation*}
a(t) = \exp (t \mathbf h) \qquad \text{with } \mathbf h = \diag \left(\frac{n - 1 }{ 2}, \frac {n - 3 }{ 2}, \dots, \frac {- n + 1 }{ 2} \right) \in \mathfrak h
\end{equation*}
and $c _ n = [2 \binom {n+1 }{ 3}] ^{-1}$.
Since, by Ruelle's inequality\footnote{This is usually stated for diffeomorphisms of compact metric spaces. In the context we need, the claim is contained in \cite[Theorem 9.7]{Margulis-Tomanov}.} the Haar measure on $\Gamma \backslash G$ has maximal entropy, we have:
\begin{align*}
h _ \mu (a ({\cdot} )) & = \int h _ \nu (a ({\cdot} )) d \tau (\nu) \\
& = \sum_ {i \geq 1} a_i h _ {\nu _ i} (a ({\cdot} )) \leq  h _ {\haarG} (a ({\cdot} )) \sum_ {i \geq 1} a_i
\end{align*}
it follows that $\sum_ {i \geq 1} a_i \geq c _ n \rho$ as claimed.
\end{proof}

%

\subsection{Periodic orbits within a fixed compact set}

We now turn to proving Theorem~\ref{theorem about bounded orbits}.
This theorem states that if $n \geq 3$, $G= \PGL_n(\R)$, and $\Gamma = \PGL_n(\Z)$, for any $\varepsilon > 0$, and any compact $\Omega \subset X = \Gamma \backslash G$, the total volume of all periodic $H$ orbits completely contained in $\Omega$ with discriminant $\leq \Delta$ is $\ll _ {\varepsilon, \Omega} \Delta ^ \varepsilon$.

\begin{proof} [Proof of Theorem~\ref{theorem about bounded orbits}]
Fix $\Omega, \varepsilon$. Suppose in contradiction that there is some $C$ and an infinite sequence of $\Delta _ i \to \infty$ so that, for every $i$, there is a collection of periodic orbits
$ Y _ i = \left\{ y_{i 1} H, \dots, y_{i n _ i} H \right\}$ so that
\begin{enumerate}
\item[(i)] each $y _ {i j} H \subset \Omega$,
\item[(ii)] the discriminant $\disc ( y _ {i j} H) \leq \Delta _ i$,
\item[(iii)] $\sum_ {j = 1} ^ {n _ i} \vol (y _ {i j} H) \geq C \Delta _ i ^ \varepsilon$.
\end{enumerate}
Define for each $i$ a probability measure $\mu _ i$ as in Theorem~\ref{comprehensive equidistribution theorem}. By (i) above, all the measures $\mu _ i$
are supported on the compact set $\Omega$, and so without loss of generality we can assume that $\mu _ i$ converge weak$^*$ to some probability measure $\mu$ which would also be supported on $\Omega$.

However, $Y _ i$ satisfy all the assumptions of Theorem~\ref{comprehensive equidistribution theorem}, and it follows that $\mu$ cannot be compactly supported --- a contradiction.
\end{proof}

\subsection{Density of periodic orbits for $\R$-split division algebras}
Because (3) of Theorem~\ref{comprehensive equidistribution theorem} is automatically satisfied, for the compact quotients $X = \Gamma \backslash G$ arising from $\R$-split degree $n \geq 3 $ division algebras over $\Q$ one can get more precise information regarding density properties of periodic $H$-orbits. We recall Theorem~\ref{theorem about density in the compact case} in a slightly more explicit phrasing:

\newtheorem*{densitytheorem}{Theorem~\ref{theorem about density in the compact case}{$'$}}

\begin {densitytheorem}
Let $X = \Gamma \backslash G$ be as in (L-2) of the introduction, i.e. $\Gamma$
is a lattice in $G = \PGL_n(\R)$ associated to a  {\em division algebra} over $\Q$ (see Section \ref{CSANot} for details).
For any $i$ let $(x_{i,j})_{j = 1, \dots, N_i}$ be a finite collection of $H$-periodic points with distinct $H$-orbits such that
\begin{equation*}
\sum_ {j = 1} ^ {N _ i} \vol (x_{i,j}H) \geq C \max_ j
(\disc(x_{i,j}H))^ \rho .\end{equation*} Suppose that there is no
periodic $L$-orbit of a group $H < L < G$ (with both inclusions
proper) containing infinitely many $x _ {i,j}$. Then $\overline
{\bigcup_ {i,j} x_{i,j}H} = \Gamma \backslash G$.
\end {densitytheorem}

\begin{proof}
Suppose that the sequence of collections of $H$-periodic orbits
\begin{equation*}
Y _ i = \left\{ y _ {i 1} H, \dots, y _ {i j} H \right\} \qquad i=1, 2, \dots
\end{equation*}
forms a counterexample for some fixed $C, \rho$, i.e. this sequence
satisfies all the conditions of the above statement, but $\bigcup_
{i = 1 } ^ {n _ i} y _ {ij} H$ do not become dense. Then
there is
some open $U \subset X$ so that for every $i, j$ we have that $y _
{ij} H \cap U = \emptyset$.

Define probability measures $\mu _ i$ for each $Y _ i$ as in the
proof of Theorem~\ref{theorem about bounded orbits}, and passing to
a subsequence if necessary we may assume that the measures $\mu _ i$
converge in the weak$^*$ topology to a probability measure $\mu$. By
Theorem~\ref{comprehensive equidistribution theorem}, $\mu = \sum_
{i = 0}^ \infty a_i \nu _ i$ with $a_i \geq 0$, \ $h _ {\nu _ 0} (a
({\cdot})) = 0$ and each $\nu _ i$ an $L_i$ invariant probability
measure on a single periodic $L_i$ orbit, $L_i >H$ a closed subgroup
of $G$, and $\sum_ {i \geq 1} a _ i \geq c _ n \rho$.

Suppose $\nu _ 1$ is a $L _ 1$ invariant probability measure on the periodic orbit $x_1 L_1$, with $a_1 > 0$ and $L _ 1 >H$. It follows from \cite[Thm.~2.13]{Prasad-Raghunathan} (see the beginning of the proof of Theorem~\ref{Lindenstrauss-Barak}) that $x _ 1 L _ 1$ contains a periodic $H$ orbit, say $z H$.

Define
\begin{equation*}
Y = \bigcap_ {k \geq 1} \overline {\bigcup_ {i \geq k} \bigcup_ {j =
1} ^ {n _ i} y _ {ij} H} .\end{equation*} Since $U$ is open and
disjoint from all $y _ {ij} H$, we have that $U \cap Y = \emptyset$.
Also,  since $\mu _ i$ converge to $\mu$, and $\mu \geq a_1 \nu _
1$, we know that $z L_1 = x_1 L_1 \subset Y$. Let $M \geq L _ 1$ be
a maximal closed connected subgroup of $G$ so that the orbit $zM$ is
periodic and is contained in $Y$. Clearly $U \cap Y = \emptyset$
implies that $M \neq G$. On the other hand, since $M \neq G$ it
follows from the assumptions of Theorem~\ref{theorem about density
in the compact case} that there is some $i_0$ so that $y_{ij} \not
\in zM$ for all $i \geq i _ 0$, hence
\begin{equation*}
Y = \bigcap_ {k \geq i_0} \overline {\left (\bigcup_ {i \geq k} \bigcup_ {j = 1} ^ {n _ i} y _ {ij} H \right)} \subset \overline {Y \setminus zM}
.\end{equation*}
By Theorem~\ref{Lindenstrauss-Barak}, there is a strictly bigger $\tilde M > M$ so that $z \tilde M$ is periodic and contained in $Y$ --- a contradiction.
\end{proof}

\section{Applications to sharpening Minkowski's theorem} \label{minkowski}

In this Section, we translate one of the foregoing results into number-theoretic terms.
As it turns out, it has a rather pleasant application to sharpening an old result of Minkowski.

\subsection{Minkowski's theorem.}
We recall Minkowski's theorem regarding ideal classes, which in particular implies finiteness of the ideal class group:

\begin{Theorem}[Minkowski's theorem] \label{Minkowski theorem}
Let $ K$ be a number field with maximal order $\order _ K$.  Then any ideal class for $K$ possesses
a representative $J \subset \order _ K$ of norm $N(J)=O(\sqrt{\disc(K)})$ where the implicit constant depends only on $d$.
\end{Theorem}

We conjecture that this is not sharp for totally real number fields of degree $d \geq 3$
insofar as one can replace $O(\cdot)$ by $o(\cdot)$:

\begin{Conjecture} \label{conj:mink}
Suppose $d \geq 3$ is fixed. Then any ideal class in a totally real number fields of degree $d$ has a representative
of norm $o(\sqrt{ \disc(K)})$.
\end{Conjecture}
We expect that this is {\em false} for $d=2$. See discussion in Sec. \ref{sec:traduction}.

Let $K$ be a number field, $\order _ K$ its integer ring, and $[J]$ an ideal class of $\order_ K$. We will denote the regulator of $K$ (more precisely of the integer ring $\order_K$) by  $R_K$. Define
\begin{align*}
m ([J],K) & = \min_ {J ' \in [J], J' \subset \order_K} N (J)\\
m(K) & = \max_ {[J] } m ([J], K)
.\end{align*}
where, in the latter definition, the maximum is taken over all ideal classes of $K$.
Let $h_{\delta}(K)$ be the number of ideal classes in $K$ with $m ([J],K) > \delta \disc (K) ^ {1/2}$.

We prove the following towards Conjecture~\ref{conj:mink}:

\begin{Theorem} \label{thm:mink}
Let $d \geq 3$, and let $K$ denote a totally real number field of degree $d$.
For all $\varepsilon, \delta > 0$ we have:
\begin{equation} \label{eq:good}
\sum_{\disc(K) < X}
R_K h_{\delta}(K) \ll_{\varepsilon, \delta} X ^{\varepsilon}
\end{equation}
In particular:
\begin{enumerate}
\item  ``Conjecture \ref{conj:mink} is true for almost all fields'': The number of fields $K$ with discriminant $\leq X$ for which $m(K) \geq \delta \cdot \disc(K)^{1/2}$
is $O_{\epsilon}(X ^{\epsilon})$, for any $\epsilon, \delta >0$;
\item ``Conjecture \ref{conj:mink} is true for fields with large regulator'':
If $K_i$ is any sequence of fields for which $\liminf \frac{\log R_K}{\log \disc(K)} > 0$,
then $m(K_i) = o(\disc(K_i)^{1/2})$.
\end{enumerate}
\end{Theorem}

By comparison, we note that Conjecture~\ref{conj:mink} follows from Conjecture~\ref{MargulisConjecture} -- see Corollary \ref{Conjectureimplication}.  Moreover, if one assumes the GRH one may show
that $m(K) \ll_{\varepsilon} (\disc K)^{1+\varepsilon} R_K ^{-2}$, i.e., one may show a quantitative version
of the second assertion of the Theorem, but {\em only} in the range when $R_K$
is significantly larger than $(\disc K)^{1/4}$.  The proof of this is very easy, but we postpone
it to the paper \cite{ELMV3}, where we discuss $\zeta$-functions more systematically.

While the statement of both the Conjecture and Theorem are quite
modest (postulating only $o(D ^{1/2})$ instead of $O(D ^{1/2})$), we
note that Minkowski's bound is in any case very close to sharp:
\begin{Proposition} \label{Close to sharp}
For any $d \geq 2$ there exists a $c' >0$ such that there is an infinite set of totally real fields of degree $d$ for which $m(K) \geq c' \cdot \disc(K)^{1/2} (\log \disc(K))^{1-2d}$.
\end{Proposition}

\begin{proof}
It has been proved by Duke \cite{DukeCompositio} that
there exist infinitely many totally real fields $K$ of degree $d$
whose class number is $\geq c(d) \disc(K)^{1/2} (\log \disc(K))^{-d}$, where $c(d)$ is an explicit positive function of $d$. On the other hand, the total number of integral ideals with norm $m$
is bounded by $\sigma_{d}(m)$, where $\sigma_d(m)$ is the number of ways of writing $m$ as
an ordered product of $d$ non-negative integers. 
So, the number of integral ideals with norm $\leq X$ is bounded by $\sum_{m \leq X} \sigma_{d}(m)
\ll X \log(X)^{d-1}$. Thus there must exist at least one ideal class which has no representative
with norm $\ll \disc(K)^{1/2} (\log \disc(K))^{1-2d}$.
\end{proof}

\subsection{Translation to dynamics.} \label{sec:traduction}
Throughout this section, let $G=\PGL_d(\R)$,\ $\Gamma = \PGL_d(\Z)$, \ $X = \Gamma \backslash G$, and $H<G$ the group of diagonal matrices.

As we show in Proposition~\ref{Minkowski and dynamics proposition} below, if $K$ is a totally real number field of a degree $d$ and $\theta: K \otimes \R \to \R ^ d$ is an algebra isomorphism, $m ([J],K)$ is intimately related to how far the $H$-orbit of (the homothety class of) the lattice $\theta (J ^{-1})$ (considered as an element of $X$) penetrates the cusp of $X$. Recall that this orbit
is the periodic $H$-orbit associated to the data $(K, [J ^{-1}], \theta)$, in the notation of Corollary~\ref{classical bijection}.

For $d=2$ there is no reason to believe there should be any constraints on such an $H$-orbit, and we expect that Theorem~\ref{Minkowski theorem} cannot be improved in this case. Establishing this rigorously is somewhat complicated because not all periodic $H$-orbits can be obtained as $\theta (J ^{-1})$ for some ideal $J \subset \order _ K$ (by Corollary~\ref{classical bijection} the periodic $H$-orbit are encoded by triples $(K, [L], \theta)$, and in this section we are only interested in periodic orbits for which the associated order $\order_L$ is the maximal order, i.e. $\order _ K$).

Because the dynamics of the action of $H$ on $X$ is much more rigid for $d \geq 3$ (see \S\ref{dynamics section}), we expect the stronger Conjecture~\ref{conj:mink} to hold, and we show in \S\ref{sharp Minkowski proof} that Conjecture~\ref{conj:mink} would follow from Conjecture~\ref{MargulisConjecture}.

\subsection{Small norm representatives of ideal classes and $H$-orbits}
\label{sharp Minkowski proof}

For any $\delta > 0$ let $\Omega ' _ \delta \subset X$ denote the set of homothety classes of lattices $\Lambda < \R ^ d$ containing no vector $v$ with $\norm {v} _ \infty ^ d < \delta \operatorname{covol} (\Lambda)$; equivalently
\begin{equation*}
\Omega ' _ {\delta} = \left\{ \Gamma g \in X : \norm { m g }^ d _
\infty \geq  \delta \det (g) \text { for every $m \in \Z ^ d$}
\right\} .\end{equation*} This set is compact by Mahler's
Compactness Criterion \cite [Cor.~10.9]{Raghunathan}, giving a
slightly different system of neighborhoods of infinity to that used
earlier (\eqref{omegadef}).  We note that $\Omega'_{a} \subset
\Omega'_{b}$ if $a > b$. The previous system of neighbourhoods,
however, has the advantage of being defined for a general group.

\begin{Proposition} \label{Minkowski and dynamics proposition} Let $K$ be a totally real number field of degree $d$, $\theta _ K: K \otimes \R \to \R ^ d$ an algebra isomorphism, and $J$ a fractional ideal of $K$.
Let $Y$ be the periodic $H$-orbit on $\Gamma \backslash G$ associated to the
data $(K, [J ^{-1}], \theta)$.
Then the following are equivalent:
\begin{enumerate}
\item [(a)]
$m (K, [J]) < \delta \disc (K) ^{1/2} $
\item[(b)] $Y$ is not contained in $\Omega ' _ \delta $.
\end{enumerate}
\end{Proposition}

\begin{proof}
Let us first notice that, if $\Lambda \subset \R ^ d$ is a lattice, then $H. \Lambda \subset \Omega'_{\delta}$
if and only if, for all $(x_1, \dots, x_d) \in L$, we have $\prod_i |x_i| \geq \delta \mathrm{covol}(\Lambda)$.

Apply this remark to the lattice $\Lambda = \theta(J ^{-1})$; recall that $Y$ is precisely
the $H$-orbit of this lattice (or, rather, its homothety class). The covolume of $\Lambda$
is $N(J)^{-1} (\disc K)^{1/2}$, where $N(J)$ is the norm of the ideal $J$.
On the other hand the product of the coordinates of $\theta(x)$
(for $x \in K$) is the norm $N(x)$.

Therefore, condition (b) is equivalent to:
\begin{equation} \label{normcond} \mbox{ There exists }x \in J ^{-1}, \ \ \ |N(x)| \leq \delta  N(J)^{-1} (\disc K)^{1/2}.\end{equation}

But elements $x \in J ^{-1}$ map surjectively onto the set of ideals $I \subset \order_K$ which belong to the same ideal class as $J$; the map is $x \mapsto x. J$.  Therefore, \eqref{normcond} translates
to condition (a).
%
%
\end{proof}

Note that, for $\delta$ sufficiently large, $\Omega'_{\delta}$ is empty. Thus the Proposition
immediately implies Theorem \ref{Minkowski theorem} (at least for totally real $K$).
Another direct corollary of Proposition \ref{Minkowski and dynamics proposition} is the following:

\begin{Corollary} \label{Conjectureimplication}
Conjecture~\ref{MargulisConjecture} implies Conjecture~\ref{conj:mink}.
\end{Corollary}

\begin{proof}
If  Conjecture~\ref{conj:mink} were false, then for some $\delta > 0$ there is an infinite sequence of totally real fields $K_i$ (equipped with an algebra isomorphism $\theta _ i: K _ i \otimes \R \to \R ^ d$) and ideals $J _ i \subset \order _ {K _ i}$ with $m ([J _ i], K _ i) \geq \delta \sqrt {\disc (K)}$. By Proposition~\ref{Minkowski and dynamics proposition} this gives us an infinite sequence of periodic $H$-orbits all inside the fixed compact set $\Omega ' _ \delta$, in contradiction to Conjecture~\ref{MargulisConjecture}.
\end{proof}

Substituting Theorem~\ref{theorem about bounded orbits} for Conjecture~\ref{MargulisConjecture}, a similar argument gives the following:

\begin{proof} [Proof of Theorem~\ref{thm:mink}]
Suppose Theorem~\ref{thm:mink} was false. Then there will be $C, \varepsilon, \delta > 0$ and a sequence of integers $D_i \to \infty$
so that
\begin{equation} \label{many Greek delta bad fields}
\sum_{\disc(K) < D_i}
R_K h_{\delta}(K) > C D_i ^{\varepsilon} \qquad \text{for all $i$},
\end{equation}
the summation being on totally real number fields of degree $d$.
Fix for every totally real field $K$ for which $m(K) \geq \delta \sqrt {\disc (K)}$ an algebra isomorphism $\theta _ K: K \otimes \R \to \R ^ d$. Let $[J_{K,j}]$ ($j=1$,\dots, $h _ {\delta} (K)$) be the ideal classes of $K$ with $m ([J_{K,j}],K)\geq \delta \sqrt {\disc (K)}$.

Let $Y_{K,j}$ be the periodic $H$-orbit parameterized (in the language of Corollary \ref{classical bijection}) by $(K, [J_{K,j}^{-1}], \theta_K)$, i.e. the $H$-orbit of the homothety class of $\theta_K(J_{K,j}^{-1})$.  Corollary \ref{classical bijection} states that $Y_{K,j}$ are periodic $H$-orbits
with volume proportional to $R_K$ and discriminant proportional to $\disc(K)$;
Proposition \ref{Minkowski and dynamics proposition} shows that $Y_{K,j} \subset \Omega'_{\delta}$.

%
Our assumption \equ{many Greek delta bad fields} implies that the collections of periodic $H$-orbits
\begin{equation*}
\mathcal{C} _ i = \left\{ Y_{K,j}: \text{$K$ totally real with $\disc (K) < D_i$}, 1 \leq j \leq h_\delta(K) \right\}
\subset \Omega'_{\delta}
\end{equation*}
have discriminant $\leq D_i$ and total volume $\gg D_i ^{\varepsilon}$,
contradicting Theorem~\ref{theorem about bounded orbits}.
\end{proof}

\section{Examples of periodic orbits which do not equidistribute}\label{examples section}

In this section we discuss in detail the example of $n=2$ (where
rigidity for the action of the diagonal group completely breaks
down) as well as examples for $n \geq 3$ showing that individual
orbits can escape to the cusp. The latter class of examples are
particularly relevant, because they suggest strongly that some of
the hypotheses in our previous theorems cannot be easily removed.

\subsection{Abundance of periodic orbits in compact sets in $\PGL_2(\Z) \backslash \PGL_2(\R)$}\label{rank one section}

As we discussed in section \ref{sec:gn} there are two natural
parameters attached to each periodic $H$-orbit $ \Gamma g H$  of an
arithmetic quotients $\Gamma \backslash G$. The discriminant
$\disc(\Gamma g H)$ which measures the arithmetic complexity of the
orbit and the volume (or regulator) of the orbit $\vol(\Gamma g H)$
which is the covolume of the lattice $g ^{-1}\Gamma g \cap H$ in
$H$.  As we have seen it is natural for us to order the periodic
orbits by their discriminant but to use the volume as weights when
counting the orbits. In this sense the next theorem shows that quite
many periodic orbits under the geodesic flow on $\PGL_2(\Z)
\backslash \PGL_2(\R)$ belong to a fixed compact set.

\begin{Theorem} \label{Rank one}
Let $G=\PGL_2(\R)$, $\Gamma = \PGL_2(\Z) $, and $X_{2}= \Gamma \backslash G$. Let
\begin{equation*}
H =\left\{\begin{pmatrix}e ^{t/2}&\\&e ^{-t/2}\end{pmatrix}:t \in \R\right\}
\end{equation*}
be the diagonal Cartan so that the action of $H$ on $X_2$ is
precisely the geodesic flow on the (unit tangent bundle of the)
unimodular surface. Then for every $\epsilon > 0$ there exists a
$\delta > 0$ such that
\begin{equation*}
\sum_ {\disc(\Gamma g H) \leq \Delta, \Gamma g H \subset \Omega'_{\delta}}
\vol(\Gamma g H) \gg \Delta ^{1 - \epsilon},
\end{equation*}
where $\Omega_{\delta}'$ is the compact set defined in Section
\ref{sharp Minkowski proof}.
\end{Theorem}
In comparison, we remark that without the restriction on the compact
set (cf. \cite{Siegel}) 
\begin{equation*}
\sum_ {\disc(\Gamma g H) \leq \Delta} \vol(\Gamma g H)
\ll_{\varepsilon} \Delta ^{3/2+\varepsilon}.
\end{equation*}

This theorem is just a special case of the general philosophy that
the dynamics of $\R$-rank one Cartan groups $H \subset \PGL_2$ is
very flexible.  Using similar methods the theorem can be generalized
to periodic orbits arising e.g.~from cubic number fields with one
real and one complex embedding (where $\mathbb{S} ^ 1 \times \R ^ {>
0} \cong H < \PGL_3(\R)$). We restrict our attention to the special
case of the geodesic flow only for brevity, and to be able to give a
simple self contained treatment. We will indicate references for the
general case. However, as we have seen the higher rank case is very
different, see Section \ref{applications section}.

For the proof of Theorem \ref{Rank one} we will need several lemmata.

\begin{Lemma}  \label{Regulator-discriminant-inequality}
(cf. Proposition \ref{proposition about regulator}) There exists a constant $C > 0$ such that
\begin{equation*}
\vol(\Gamma g H) \geq \log \disc (\Gamma g H) - C
\end{equation*}
for any periodic orbits $\Gamma g H$.
\end{Lemma}

\begin{proof}
Since both the volume (i.e. regulator) and the discriminant only depend on the order
and not on the proper ideal associated to the orbit in Corollary \ref{classical bijection}, we can assume without loss of generality (replacing, in the notation of that Corollary, ``$L$'' by ``$\order_L$'') that the orbit is defined by data $(\Q(\sqrt{D}), \order, \theta)$, 
where $\order \subset \Q(\sqrt{D})$ is the order of discriminant $D$.

By definition $\order$ gives rise to a lattice $ \Lambda _ \order := \theta(\order) \subset \R ^ 2 $ of covolume $ D ^{1/2}$. Since $1 \in \order$ we have $(1, 1) \in \Lambda _ \order$ which is a ``short vector''
in comparison to the covolume.  To see better why $(1, 1)$ is short, let us renormalize $ \Lambda _ \order$ by a homothety to covolume one.  Then $(D ^{-1/4},D ^{- 1/4})$
belongs to the so obtained  lattice $ \Lambda _ \order ^ 1 $.  Applying the geodesic flow
$h_t=\begin{pmatrix}e ^{t/2}&\\&e ^{-t/2}\end{pmatrix}$ (in either direction) makes this vector longer.  However, as long as $|t/2| \leq \log(D ^{1/4}/2)$ the vector
$v _ t = (D ^{- 1 / 4},D ^{- 1 / 4}) h_t ^{-1}$ would still be of norm less than one.  Note that $v _ t$ is moving along a hyperbola for varying values of $t$.
The unimodular lattice $ \Lambda _ \order ^ 1 h _ t ^{-1}$ can only contain (up to sign) one vector of length less than one.  Therefore,
$t \mapsto \Lambda_\order ^ 1 h _ t ^{-1}$ is injective for $t \in[-\frac{1}{2} \log D + 2 \log 2,\frac{1}{2} \log D - 2  \log 2]$ and the lemma follows.
\end{proof}

We recall the notion of topological entropy for a
continuous transformation $T:Y \rightarrow Y$ on a compact
metric space $(Y,d)$. For $\eta>0$ and a positive integer $N$
let $s_{\eta,N}(\beta)$  be the maximal number of points
$y_1,y_2,\ldots, y_{s_{\eta,N}}$ such that for any $i \neq j$ there exists some
$0 \leq k<N$ with $d(T ^ ky_i,T ^ ky_j)>\eta$. Then the
{\em topological entropy} is defined by
\[
\htop(T)=\lim_{\eta \rightarrow 0}
\limsup_{n \rightarrow \infty}\frac{\log s_{\eta,n}(\beta)}{n}.
\]

The following is a simple special case of a very general fact, see \cite{Kleinbock-Margulis} and
the references therein.

\begin{Lemma} \label{rank-one-entropy}
For any $ \epsilon > 0$ there exists $ \delta > 0$ such that the restriction $T$
of (the time one map  of the geodesic flow) $h _ 1$ to
\begin{equation*}
Y _ \delta = \{ x \in  \Omega'_\delta: h _ t . x \in \Omega'_\delta\mbox{ for all } t>0 \}
\end{equation*}
satisfies $\htop(T)>1-\epsilon$. In fact, we have $s_{\eta, N} \geq e ^{(1-\epsilon)N}$ for sufficiently small $\eta$.
\end{Lemma}

\begin{proof}
Recall that $u \in \R$ is badly approximable ($BA$)
if\footnote{In this context, $\langle x \rangle$ denotes the nearest integer.}  $\liminf _{ m \rightarrow \infty} m \langle m u \rangle > 0$ and that
the set of badly approximable numbers has Hausdorff dimension one. Moreover, since
the Hausdorff dimension of a countable union is a supremum of the dimensions of the individual sets, it follows that for small enough $\delta > 0$ the set $BA(\delta )$ of elements $u \in [0, \frac{1}{2}]$ with $\liminf _{ m \rightarrow \infty} m \langle m u \rangle > \delta $ has Hausdorff dimension greater than $1 - \epsilon$.  In particular, for sufficiently small $\eta _ 1 > 0$ there exist at least
$\eta _ 1 ^{- (1 - \epsilon)}$ many such $u \in BA(\delta )$ that are at least $\eta _ 1$ far apart.

Let $x _ u = \Gamma\begin{pmatrix} 1 & -u\\&1\end{pmatrix}$ for some $u \in BA(\delta )$.  We claim that $x _u \in Y_\delta$.  For, if not, there exists $t > 0$
is sothat the lattice $\Z ^ 2\begin{pmatrix} 1 & -u\\&1\end{pmatrix}\begin{pmatrix} e ^{-t/2} & \\ &e ^{t/2}\end{pmatrix}$ contains a vector
$(me ^{-t/2},| m u + m '|e ^{t/2})$ whose coordinates have product $\leq \delta$ in absolute value. Then $m \ne 0$ and
$|m( m u + m ') | < \delta $.

Since $h _ 1\begin{pmatrix} 1 & u\\&1\end{pmatrix} h _ 1 ^{- 1} = \begin{pmatrix} 1 & eu\\&1\end{pmatrix}$ it follows that the distance between two close by
points on the same orbit of $\begin{pmatrix} 1 & \R \\ &1\end{pmatrix}$
is uniformly increasing until the distance is bigger than the injectivity radius of the
local isomorphism between $G$ and $X$ at the image point. For $u _ 1, u _ 2 \in BA(\delta )$ the points $x _ i = \Gamma\begin{pmatrix} 1 & -u _ i \\ &1\end{pmatrix} $
will stay forever in $\Omega'_\delta$.  So if $u _ 1, u _ 2$ are $\eta_1 = e ^{- N} \eta _ 0$ far apart, we can ensure by choosing $\eta _ 0$ sufficiently small that $d(T ^ k.  x _ 1,T ^ k.  x _ 2) > \eta$ for some $k$ with $0 \leq k < N$.  This shows that $s _{\eta, N} \geq e ^{(1 - \epsilon) N}$ for sufficiently large $N$ and so $\htop (T) \geq 1 - \epsilon$ as claimed.
\end{proof}

To obtain periodic points from the above lemma regarding topological entropy we need the following standard facts from hyperbolic dynamics.  We do not give here the statements
in their full strength --- not even for the special case considered here, see e.g.~ \cite [Sec.~18.1.]{Katok-book}.

\begin{Lemma}  \label{Shadowing}
(Shadowing lemma) For any $\eta _ s> 0$ there exists $ \rho _ s> 0$ such that
if $d (x _ -, x _ +) < \rho_s$, then there exists $y$ with $d (h _ t. y, h _ t.x _ -) < \eta_s$ for all $ t \leq 0$ and $d (h _ t. y, h _ t.x _ +) < \eta_s$ for all $ t \geq 0$.

(Anosov closing lemma) For any $\eta _ c> 0$ there exists $ \rho _ c> 0$ such that
if $d (x, h _ N x) < \rho_c$ for some $N \geq 1$, then there exists $y$ and $T \in [N - \eta_c, N + \eta_c]$ such that $ h _ T. y = y$ and $d (h _ t. x, h _ t. y) < \eta_c$ for $t \in [0, N]$.
\end{Lemma}

\begin{proof}
Write $x _ - = g . x _ +$ with $g = \begin{pmatrix} g_{11} & g_{12}\\g_{21}&g_{22}\end{pmatrix}$, $d (g,e) < \rho$ and $ \det g = 1$.  Define
$ y = \begin{pmatrix} 1 & u\\ &1\end{pmatrix}.x _ -$.  It is easy to see that for small enough values of $u$ this $y$
will always satisfy the first part of the shadowing statement. Since $y = \begin{pmatrix} 1 & u\\ &1\end{pmatrix}
\begin{pmatrix} g_{11} & g_{12}\\ g_{21}& g_{22}\end{pmatrix}.x _ + =\begin{pmatrix} g_{11}+g_{21} u& g_{12}+g_{22} u\\ g_{21}&g_{22}\end{pmatrix}.x _ +$,
it follows similarly that for $u = - \frac {g_{12}} {g_{22}}$ both parts to shadowing hold true.

Suppose now $h _ N .x = g .x$ with $g$ as before.  Define $ \tilde x = \begin{pmatrix} 1 & \\ u&1\end{pmatrix}.x$.  Then
\begin{align*}
h _ N .\tilde x & = h _ N\begin{pmatrix} 1 & \\ u&1\end{pmatrix}h _ N ^{- 1}h _ N. x =
\begin{pmatrix} 1 & \\ e ^{-N}u&1\end{pmatrix} g\begin{pmatrix} 1 & \\ -u&1\end{pmatrix}. \tilde x \\
& =
\begin{pmatrix} g_{11}-ug_{12} & g_{12}\\ -e ^{-N}g_{12}u ^ 2+e ^{-N}g_{11}u-g_{22}u+g_{21}&g_{22}+e ^{-N}ug_{12}\end{pmatrix}.\tilde x.
\end{align*}
If $ \rho_c$ and so $(g_{11}-1),g_{12},g_{21},(g_{22}-1)$ are sufficiently small, the quadratic polynomial in the low left corner
has a unique root $u$ close to zero. In other words by replacing $x$ by the close by point $\tilde x$ we can assume that $g_{21} = 0$.
Now define $y = \begin{pmatrix} 1 & v\\ &1\end{pmatrix}. \tilde x$. A similar calculation (which in fact now leads to a linear equation) shows the Anosov closing lemma.
\end{proof}

\begin{Lemma} \label{many periodic orbits}
Let $\epsilon$ and $\delta$ be as in Lemma \ref{rank-one-entropy}.
Then the number of periodic orbits of length less than $V$ that are contained
in $\Omega'_\delta$ is $\gg e ^{(1-\epsilon)V}$.
\end{Lemma}

\begin{proof}
By Lemma \ref{rank-one-entropy} we can choose $\eta$ such that for large enough $N$ there are $\gg e ^{(1-\epsilon)V}$
many points $x_1,\ldots \in Y_\delta$ with $d(h_k.x_i,h_k.x_j)>\eta$ for $0 \leq k=k(i,j) < N$.
We are going to apply both parts of Lemma \ref{Shadowing}. The parameters $\eta_s$ and $\eta_c$
we choose such that after application of both statements to get points $y_1,\ldots \in X$ we still have $d(h_k.y_i,h_k.y_j)>\eta/2$ for $k=k(i,j)$.
By choosing the parameters for the closing statement first, we can assume that $2 \eta_s < \rho_c$.  Since the geodesic flow is ergodic,
there exists a finite $\rho_s$-dense set of points $z_1,\ldots,z_\ell$ in $\Omega'_\delta$ and for every $z_i$ finitely many times $t(i,j)$ such that
the set of points $h_{t(i,j)}.z_i$ for varying $j$ is $(\rho_c-2 \eta_s)$-dense in $\Omega'_\delta$. Therefore, for every $x_m$ we can find a $z_i$
with $d(h_N.x_m,z_i)<\rho_s$ and obtain by shadowing a point $\tilde x_m$ with $d(\tilde x_m, h_{N+t(i,j)}.\tilde x_m)<\rho_c$.
Applying the closing lemma, we arrive at points $y_1,\ldots \in X$ with periodic orbits. These periodic points belong for varying $N$ to
a fixed compact set, namely to a neighborhood of the union of $\Omega'_\delta$ and the finite pieces of the orbits of $z_1,\ldots,z_\ell$ used in the argument.
The length of the orbits are all below $V=N+\max_{i,j}t_{i,j}+\eta _c$. In particular, every periodic orbit can be counted at most $O(V)$ times since any
two points in the list are at some moment $\eta/2$ far apart. Since $\epsilon$ is arbitrary, this does not affect the statement of the lemma.
\end{proof}

\begin{proof}[Proof of Theorem \ref{Rank one}]
Let $\epsilon, \delta$ be as in Lemma \ref{rank-one-entropy}. Since there is a lower bound
on the length of an $H$-orbit we have
\begin{equation*}
\sum_ {\vol(\Gamma g H) \leq V, \Gamma g H \subset \Omega'_\delta}
\vol(\Gamma g H) \gg e ^{(1-\epsilon)V}
\end{equation*}
by Lemma \ref{many periodic orbits}.
If we set $V=\log \Delta -C$, then $\vol(\Gamma g H) \leq V$ implies $\disc(\Gamma gH)\leq \Delta$ by Lemma \ref{Regulator-discriminant-inequality}. Therefore,
\begin{equation*}
\sum_ {\disc(\Gamma g H) \leq \Delta, \Gamma g H \subset \Omega'_\delta}
\vol(\Gamma g H) \gg e ^{(1-\epsilon)V} \gg \Delta ^{1-\epsilon}
\end{equation*}
as claimed.
\end{proof}

\subsection {Fields of small regulator and non-equidistributed orbits in higher rank; proof of Theorem
\ref{badcompactorbits}.}  \label{Dukefields}

As an example which is illuminating in its own right, we construct using the above parameterizations
a sequence of periodic $H$-orbits for which a positive fraction of mass escapes to the cusp. More precisely, we shall construct
a sequence of periodic $H$-orbits $Y_i$ on $\PGL_n(\Z) \backslash \PGL_n(\R)$ with the following properties, for some constants $a_n, b_n >0$:
\begin{enumerate}
\item The number of $Y_i$ with discriminant $\leq \Delta$ is $\gg \Delta ^{a_n}$.
\item Any weak limit of the $H$-invariant  probability measures $\mu_{Y_i}$ on $Y_i$
has total mass $\leq 1-b_n$.
\end{enumerate}

This example will, in particular, prove Theorem \ref{badcompactorbits} asserted in the Introduction.
We would like to thank Peter Sarnak for suggesting the possibility of it.

%

\begin{proof}
Indeed, Duke has established (cf. \cite[\S 3, Prop. 1] {DukeCompositio}) that, for every degree $n$, there
is a constant $C = C(n)$ and infinitely many totally real fields
$K_i$ whose Galois closure has Galois group $S_n$ and whose regulators satisfy the inequality $R_{K_i} < C (\log \disc (K_i))^{n-1}$.

For each such $K_i$, let $\order_i$ be the maximal order, fix
an algebra isomorphism $K_i \otimes \R \rightarrow \R ^ n$, and consider
the periodic orbit associated to the triple $(K_i, [\order_i], \theta)$.

Let $\alpha: K_i \otimes \R \rightarrow \R ^ n$ be the map
$x \mapsto (\log |\theta(x)_1|, \log |\theta(x)_2|, \dots, \log |\theta(x)_n|)$.
Let $\Lambda$ be the image of $\order_i ^{\times}$ under the map $\alpha$.
Then the $\R$-span $\Lambda_{\R}$ of $\Lambda$ is, by Dirichlet's unit theorem, precisely the subspace
$\{(y_1, \dots, y_n) \in \R ^ n: \sum y_i = 0 \}$, the image under $\alpha$
of elements $(K_i \otimes \R)^{(1)}$ of norm $1$.  Moreover, $\alpha$ induces an isomorphism of torii
$$\alpha: (K_i \otimes \R)^{(1)}/\order_i ^{\times} \rightarrow \Lambda_{\R}/\Lambda$$

Endowing the left-hand torus with the measure transported from Lebesgue measure on $\Lambda_{\R} \subset \R ^ n$, the volume of the left-hand torus is $R_{K_i}$. On the other hand:
\begin{enumerate}
\item
Since $K_i$ has no intermediate subfields, all nonzero elements of
$\Lambda$ have Euclidean length $\geq c \log \disc(K_i)$, for a
certain constant $c$ depending only on $n$ (cf. e.g. \cite{Remak} or
proof of Corollary \ref{Cor:Explication} (ii))
\item If $x \in (K_i \otimes \R)^{(1)}$ is such that $|\log |\theta(x)_i|| \leq \frac{c}{2} \log \disc(K_i)$ for each $i$,
then the lattice $\theta(x . \order_i)$, which has covolume $\asymp \disc(K)^{1/2}$, contains a vector of Euclidean length
$\leq \sqrt{n} (\disc K_i)^{c/2}$.
\end{enumerate}

Fix a compact subset $\Omega \subset \PGL_n(\Z) \backslash
\PGL_n(\R)$. It follows from these remarks that there is a subset
$(K_i \otimes \R)^{(1)}/\order_i ^{\times}$ of volume $\gg c' (\log
\disc (K_i))^{n-1}$, where $c'$ depends only on $n$, such that the
lattice $\theta(x . \order_i)$ does not belong to $\Omega$ for all
sufficiently large $i$.   Since the total volume of $(K_i \otimes
\R)^{(1)}/\order_i ^{\times}$ w.r.t. the measure defined above was
$R_{K_i} < C (\log \disc(K_i))^{n-1}$, this implies the conclusion.
\end{proof}

\bibliographystyle{plain}
\bibliography{ELMV1-final}

\end{document}